\documentclass[12pt,leqno,a4paper]{article}

\usepackage{amsmath}
\usepackage{amssymb}
\usepackage{amsthm}
\usepackage[UKenglish]{babel}
\usepackage{bbm}               % for \tfn
\usepackage{ifthen}
\usepackage{sectsty}
\usepackage{url}

\sectionfont{\large}
\subsectionfont{\normalsize}

\mathchardef\ordinarycolon\mathcode`\:
\mathcode`\:=\string"8000
\def\vcentcolon{\mathrel{\mathop\ordinarycolon}}
\begingroup \catcode`\:=\active
  \lowercase{\endgroup
  \let :\vcentcolon
  }

\theoremstyle{plain}
\newtheorem{thm}{Theorem}[section]

\newtheorem{lem}[thm]{Lemma}
\newtheorem{prop}[thm]{Proposition}

\theoremstyle{definition}
\newtheorem{defn}[thm]{Definition}
\newtheorem{exmp}[thm]{Example}
\newtheorem{notation}[thm]{Notation}
\newtheorem{rem}[thm]{Remark}

\newcommand{\agb}{\mathcal{M}}
\newcommand{\bebe}{\Gamma}
\newcommand{\bbebe}{{\widetilde{\bebe}}}
\newcommand{\csa}{\mathcal{A}}
\newcommand{\eevecs}{{\widetilde{\evecs}}}
\newcommand{\elltwo}{{L^2(\R_+;\mul)}}
\newcommand{\eps}{\epsilon}
\newcommand{\evec}[1]{\evecc(#1)}
\newcommand{\evecc}{\varepsilon}
\newcommand{\evecs}{\mathcal{E}}
\newcommand{\expn}{\mathbb{E}}
\newcommand{\ffock}{{\widetilde{\fock}}}
\newcommand{\fock}{\mathcal{F}}
\newcommand{\half}{\mbox{$\frac12$}}
\newcommand{\hilb}{\mathsf{H}}
\newcommand{\hilc}{\mathsf{K}}
\newcommand{\id}{I}
\newcommand{\indf}[1]{\indff_{#1}}
\newcommand{\indff}{1}
\newcommand{\ini}{\mathsf{h}}
\newcommand{\lind}{L}
\newcommand{\mmodf}[1]{\overline{\overline{#1}}}
\newcommand{\mmul}{{\widehat{\mul}}}
\newcommand{\modf}[1]{\overline{#1}}
\newcommand{\mul}{\mathsf{k}}
\newcommand{\opsp}{\mathsf{V}}
\newcommand{\scale}{s}
\newcommand{\uhm}{u}
\newcommand{\Vac}{\Omega}
\newcommand{\vac}{\omega}
\newcommand{\vhm}{v}
\newcommand{\vint}{\Lambda_\Vac}

\renewcommand{\geq}{\geqslant}
\renewcommand{\leq}{\leqslant}
\renewcommand{\Pr}{\mathbb{P}}

\newcommand{\bp}{\mathbf{p}}
\newcommand{\bq}{\mathbf{q}}
\newcommand{\bt}{\mathbf{t}}

\newcommand{\C}{\mathbb{C}}
\newcommand{\R}{\mathbb{R}}
\newcommand{\Z}{\mathbb{Z}}

\newcommand{\algten}{\odot}
\newcommand{\bop}[2]%
{\ifthenelse{\equal{#2}{}}{\bopp(#1)}{\bopp(#1;#2)}}
\newcommand{\bopp}{\mathcal{B}}
\newcommand{\cb}{\mathrm{cb}}
\newcommand{\cbo}[2]
{\ifthenelse{\equal{#2}{}}{\cboo(#1)}{\cboo(#1;#2)}}
\newcommand{\cboo}{\mathcal{CB}}
\newcommand{\Cf}{\textit{Cf.~}}
\newcommand{\cf}{\textit{cf.~}}
\newcommand{\dyad}[2]{|#1\rangle\langle#2|}
\newcommand{\e}{\mathrm{e}}
\newcommand{\etc}{\textit{et cetera}}
\newcommand{\gap}{\vspace{1ex}\noindent}
\newcommand{\I}{\mathrm{i}}

\newcommand{\ie}{\textit{i.e., }}
\newcommand{\kbo}[3]
{\ifthenelse{\equal{#3}{}}{#1\kboo(#2)}{#1\kboo(#2;#3)}}
\newcommand{\kboo}{\mathcal{B}}
\newcommand{\lift}[2]{{#1}\matten{\id_{\bop{#2}{}}}}
\newcommand{\lin}{\mathop{\mathrm{lin}}}
\newcommand{\lop}[2]%
{\ifthenelse{\equal{#2}{}}{\lopp(#1)}{\lopp(#1;#2)}}
\newcommand{\lopp}{\mathcal{L}}
\newcommand{\mat}[2]{{#2}\matten\bopp(#1)}
\newcommand{\matt}{\mathrm{M}}
\newcommand{\matten}{\mathop{\otimes_\matt}}
\newcommand{\mvec}[2]%
{\left(\begin{smallmatrix}#1\\#2\end{smallmatrix}\right)}
\newcommand{\nnabla}{\widehat{\nabla}}
\newcommand{\rd}{\mathrm{d}}
\newcommand{\splten}{\otimes}
\newcommand{\std}{\,\rd}
\newcommand{\tfn}[1]{\mathbbm{1}_{#1}}
\newcommand{\uwkten}{\mathop{\overline{\otimes}}}

\begin{document}

\begin{center}
{\LARGE Random-walk approximation\\[0.25ex]
        to vacuum cocycles}\\[1.5ex]
{\large Alexander C.~R.~Belton}\\[0.5ex]
{\small Department of Mathematics and Statistics\\
Lancaster University, United Kingdom\\[0.5ex]
\textsf{a.belton@lancaster.ac.uk} \qquad \today}
\end{center}

\begin{abstract}
\noindent
Quantum random walks are constructed on operator spaces with the aid
of matrix-space lifting, a type of ampliation intermediate between
those provided by spatial and ultraweak tensor products. Using a form
of Wiener--It\^o decomposition, a Donsker-type theorem is proved,
showing that these walks, after suitable scaling, converge in a strong
sense to vacuum cocycles: these are vacuum-adapted processes which are
Feller cocycles in the sense of Lindsay and Wills. This is employed to
give a new proof of the existence of $*$-homomorphic quantum
stochastic dilations for completely positive contraction semigroups on
von~Neumann algebras and separable unital $C^*$~algebras. The
analogous approximation result is also established within the standard
quantum stochastic framework, using the link between the two types of
adaptedness.
\end{abstract}

\gap
{\footnotesize\textit{Key words:} quantum random walk; repeated
interactions; Donsker's invariance principle; functional central-limit
theorem; quantum stochastic dilation; quantum dynamical semigroup;
completely positive contraction semigroup; toy Fock space; discrete
approximation; Feller cocycle; vacuum adapted.}

\gap
{\footnotesize\textit{MSC 2000:} %
81S25 (primary);    %% Quantum stochastic calculus
46L07,              %% Operator spaces and completely bounded maps
46L53,              %% Noncommutative probability and statistics
46N50,              %% Applications in quantum physics
47D06,              %% One-parameter semigroups and linear evolution
                    %% equations
60F17 (secondary).} %% Functional limit theorems; invariance
                    %% principles

\section{Introduction}

Let $(x_k)_{k\geq1}$ be a sequence of independent, identically
distributed random variables, with zero mean and unit variance, and
let $S_m:=\sum_{k=1}^m x_k$ for all~$m\geq0$, so that
$S=(S_m)_{m\geq0}$ is the associated random walk (with $S_0=0$). If
\[
X^{(n)}_t:=\frac{1}{\sqrt{n}}(S_m+(n t-m)x_{m+1})\qquad%
\forall\,t\in{[m/n,(m+1)/n[}
\]
then Donsker's invariance principle \cite[{\S}I.8]{RoW00} implies that
the process $X^{(n)}$ (given by scaling and linear interpolation
between the points of $S$) converges in distribution, as $n\to\infty$,
to a classical Wiener process.

This result has two parts: first, the random walk is used to construct
a continuous-time process with continuous paths (or c\`agl\`ad paths,
in the original formulation \cite{Don51}); this process, suitably
scaled, is then shown to converge to Brownian motion. Below a similar
result is established for quantum random walks: the first stage
corresponds to embedding such a walk (which lives naturally on toy
Fock space) as a continuous-time process on Boson Fock space; the
second to showing that, subject to suitable scaling, this process
converges (strongly on the exponential domain) to a type of process
called a vacuum cocycle. Our work is a considerable advance on earlier
efforts, establishing stronger results with weaker hypotheses and
providing a much clearer understanding of the underlying
structure. (Cockroft, Gudder and Hudson found \cite{CGH77} another
quantum-mechanical generalisation of Donsker's invariance principle;
in their work, classical random variables are replaced by pairs of
self-adjoint operators which satisfy the canonical commutation
relations in their Weyl form.)

The framework for our investigation is the matrix-space formulation of
quantum stochastics pioneered by Lindsay and Wills \cite{LiW01}.
Section~\ref{sec:mat} contains an introduction to these ideas and
gives the construction of a quantum random walk. It generalises the
repeated-interactions approach to open quantum systems described by
Brun \cite{Bru02}, Gough \cite{Gou04} and Attal and
Pautrat \cite{AtP06}; their articles explain the physical
interpretation of this mathematical object and the reader is
encouraged to consult them.

In quantum probability, the use of toy Fock space (\ie the tensor
product of countably many copies of $\mmul:=\C\oplus\mul$, where
$\mul$ is a complex Hilbert space and the product is taken with
respect to the unit vector $(1,0)\in\mmul$) goes back to work of
Accardi and Bach, Journ\'e and Meyer, K\"ummerer, Lindsay and
Parthasarathy \textit{et alia}, with recent contributions from Attal,
Pautrat, Gough, Sinha, Sahu, Franz and Skalski, among others; see
\cite[Introduction]{Blt07} for more information, including references
to their work beyond those given herein. In Section~\ref{sec:toy} the
necessary results on this toy space and its relationship with Boson
Fock space are summarised, together with some facts about discrete
quantum stochastic integrals (\cf \cite{Blt07}). No restriction is put
on the multiplicity space $\mul$; no advantage would accrue from
insisting upon its separability.

If $\phi$ is a linear map from the concrete operator space $\opsp$
to $\mat{\mmul}{\opsp}$, the $\mmul$-matrix space over $\opsp$ (which
lies between the spatial and ultraweak tensor products of $\opsp$ with
$\bop{\mmul}{}$), then the \emph{vacuum flow} with generator $\phi$ is
a vacuum-adapted mapping process $j$ which satisfies the Evans--Hudson
equation
\begin{equation}\label{eqn:vccdef}
\langle u\evec{f},j_t(a) v\evec{g}\rangle=%
\langle u,av\rangle+\int_0^t\langle u\evec{f},%
j_s\bigl(E^{\widehat{f(s)}}\phi(a)E_{\widehat{g(s)}}\bigr)%
v\evec{g}\rangle\std s.
\end{equation}
This generalisation of the Hudson--Parthasarathy evolution equation
\[
\rd Y_t=\rd \Lambda_F(t)Y_t\biggl(=\sum\nolimits_{\alpha,\beta}%
(F^\alpha_\beta\otimes\id_\fock)Y_t\std\Lambda^\beta_\alpha(t)\biggr)
\]
is introduced in Section~\ref{sec:coc}. The vacuum flows of interest
here are strongly regular, so uniquely specified by
(\ref{eqn:vccdef}), and are Feller cocycles in the Lindsay--Wills
sense \cite[p.284]{LiW00}.

Two approximation results are presented in Section~\ref{sec:wlk}. The
primary one, Theorem~\ref{thm:main}, gives stronger conclusions under
weaker and more natural conditions than those of Attal and Pautrat
\cite[Theorem~13 \textit{et seq.}]{AtP06}. Lindsay and Wills
discovered hypotheses (extending those of Evans and of Mohari and
Sinha) which ensure cocycle multiplicativity \cite{LiW03a,LiW03b}; it
is pleasing, but not surprising, to see them appear naturally here.

A type of Wiener--It\^o decomposition is used to establish these
results: the quantum random walk is written as a sum of iterated
`discrete integrals' and each of these is shown to converge to a
multiple quantum stochastic integral, the sum of which is equal to the
limit flow. This approach gives a much better understanding (compared
to previous work) of the structure underlying the convergence of
quantum random walks to stochastic cocycles. Although the presentation
would be a little smoother if only infinite multiplicity were to be
considered, the importance of finite dimensions (particularly for
physical applications) makes worthwhile the extra effort needed to
cover this situation too.

From its very inception, quantum stochastic calculus has been utilised
to produce dilations of quantum dynamical
semigroups. Section~\ref{sec:exa} contains a new proof of the
existence of $*$-homomorphic cocycles dilating completely positive
contraction semigroups on von~Neumann algebras \cite{GoS99,GLSW03} and
unital separable $C^*$~algebras \cite{GPS00,LiW03a,LiW03b}. This was a
primary motivation for obtaining our results, which improve on those
previously obtained by Sinha \cite{Sin06} (who worked with
one-dimensional noise) and Sahu \cite{Sah08} (who worked on
von~Neumann algebras).

There are overlapping reasons, practical and philosophical, why
\emph{vacuum} adaptedness is used below, rather than the more standard
form. (In a way, the second provides an explanation for the first.) As
is well known by now, working in a vacuum-adapted set-up frequently
leads to simpler analysis; for example, quantum stochastic integrals
may be defined directly using the classical It\^o integral, rather
than the Hitsuda--Skorokhod one. Furthermore, a non-trivial vector
process which is adapted in the usual sense cannot correspond to a
process on toy Fock space; the projection from Boson Fock space into
toy Fock space is (essentially) a conditional expectation, which
averages over the intervals into which time is partitioned. In the
vacuum-adapted situation, however, the future part of a vector process
corresponds to the vacuum vector, which is invariant under such
averaging. Section~\ref{sec:ssu} contains the analogous approximation
theorem for the standard set-up, Theorem~\ref{thm:smain}, which is
deduced from its vacuum-adapted version; it is always possible to
switch freely between vacuum-adapted and standard (identity-adapted)
perspectives \cite{Blt04}.

\subsection{Conventions and notation}
All vector spaces have complex scalar field and all inner products are
linear in the second variable. The conventions and notation of
\cite{Blt07} (which follows \cite{Lin05} for the most part) are
adopted and, as far as possible, work proceeds in a coordinate-free
manner.

The algebraic tensor product is denoted by $\algten$, the usual tensor
product of Hilbert spaces and bounded operators is denoted by
$\splten$, as is the spatial tensor product of operator spaces, and
the ultraweak tensor product is denoted by~$\uwkten$.

The vector space of linear operators between vector spaces $V$ and $W$
is denoted by $\lop{V}{W}$, or $\lop{V}{}$ if $V$ equals $W$; the
identity operator on $V$ is denoted by $\id_V$. The Banach space of
bounded operators between Banach spaces $X$ and $Y$ is denoted by
$\bop{X}{Y}$, or $\bop{X}{}$ if $X$ equals $Y$. The double commutant
of a set $A\subseteq\bop{\hilb}{}$, where $\hilb$ is a Hilbert space,
is denoted by $A''$.

An empty sum or product is equal to the appropriate additive or
multiplicative unit respectively. The restriction of a function $f$
to a set $A$ (contained in the domain of $f$) is denoted by
$f|_A$. The indicator function of a set $A$ is denoted by
$\indf{A}$. Given a proposition $P$, the expression $\tfn{P}$ has the
value $1$ if $P$ is true and $0$ if $P$ is false. The sets of
non-negative integers and non-negative real numbers are denoted by
$\Z_+$ and $\R_+$ respectively.

\section{Matrix spaces and quantum random walks}\label{sec:mat}

\begin{defn}{\cite{EfR00}}
A (concrete) \emph{operator space} $\opsp$ is a closed subspace
of $\bop{\hilb}{}$ for some Hilbert space $\hilb$. Let
$\opsp^\dagger:=\{a^*:a\in\opsp\}\subseteq\bop{\hilb}{}$ denote its
\emph{conjugate} space, where ${}^*$ is the operator adjoint.

A linear map $\phi:\opsp_1\to\opsp_2$ between operator spaces
$\opsp_1\subseteq\bop{\hilb_1}{}$ and
$\opsp_2\subseteq\bop{\hilb_2}{}$ is \emph{completely bounded} if
$\|\phi\|_\cb:=%
\sup_{n\geq1}\|\phi\algten\id_{\bop{\C^n}{}}\|<\infty$, where the
linear map
\[
\phi\algten\id_{\bop{\C^n}{}}:\opsp_1\algten\bop{\C^n}{}\to%
\opsp_2\algten\bop{\C^n}{};\ a\otimes b\mapsto\phi(a)\otimes b
\]
and $\opsp_i\algten\bop{\C^n}{}$ is regarded as a closed subspace
of $\bop{\hilb_i\otimes\C^n}{}$ for $i=1$,~$2$.
\end{defn}

\begin{defn}{\cite{LiW01}}
Given an operator space $\opsp\subseteq\bop{\hilb}{}$ and a Hilbert
space $\hilc$, the \emph{matrix space}
\[
\mat{\hilc}{\opsp}:=\{T\in\bop{\hilb\otimes\hilc}{}:%
E^x T E_y\in\opsp\ \forall\,x,y\in\hilc\}
\]
is also an operator space, where
$E^x\in\bop{\hilb\otimes\hilc}{\hilb}$ is the adjoint of the map
$E_x:y\mapsto y\otimes x$. Note the inclusions
$\opsp\splten\bop{\hilc}{} \subseteq \mat{\hilc}{\opsp} \subseteq %
\opsp \uwkten \bop{\hilc}{}$, with the latter an equality if $\opsp$
is ultraweakly closed. Note also that
$\bigl(\mat{\hilc}{\opsp}\bigr)^\dagger=\mat{\hilc}{\opsp^\dagger}$.
\end{defn}

\begin{prop}\label{prp:hsineq}
If $\hilb_1$, $\hilb_2$ and $\hilc$ are Hilbert spaces,
$T\in\bop{\hilb_1\otimes\hilc}{\hilb_2\otimes\hilc}$ and
$\{e_i:i\in I\}$ is an orthonormal basis for $\hilc$ then
\begin{equation}\label{eqn:hsineq}
\|T\|^2\leq\sum_{i,j\in I}\|E^{e_i}T E_{e_j}\|^2.
\end{equation}
Consequently, if $\phi\in\bop{\opsp_1}{\opsp_2}$, where $\opsp_1$ and
$\opsp_2$ are operator spaces, and $L$ is a finite-dimensional Hilbert
space then $\|\phi\algten\id_{\bop{L}{}}\|\leq(\dim L)\|\phi\|$.
\end{prop}

\begin{defn}
Let $\hilc\neq\{0\}$ be a Hilbert space. A linear map $\phi$ between
operator spaces is \emph{$\hilc$ bounded} if $\|\phi\|_\hilc<\infty$,
where
\[
\|\phi\|_\hilc:=\left\{\begin{array}{ll}
(\dim\hilc)\|\phi\|&\mbox{if }\dim\hilc<\infty,\\[1ex]
\|\phi\|_\cb&\mbox{if }\dim\hilc=\infty.
\end{array}\right.
\]
Let $\kbo{\hilc}{\opsp_1}{\opsp_2}$ (or $\kbo{\hilc}{\opsp_1}{}$, if
$\opsp_1$ equals $\opsp_2$) denote the collection of $\hilc$-bounded
operators between the operator spaces $\opsp_1$ and $\opsp_2$.
\end{defn}

\begin{thm}\label{thm:lift}
Let $\opsp_1\subseteq\bop{\hilb_1}{}$ and
$\opsp_2\subseteq\bop{\hilb_2}{}$ be operator spaces and let $\hilc$
be a Hilbert space. If $\phi\in\kbo{\hilc}{\opsp_1}{\opsp_2}$ then
there exists a unique map
$\lift{\phi}{\hilc}:\mat{\hilc}{\opsp_1}\to\mat{\hilc}{\opsp_2}$ such
that
\begin{equation}\label{eqn:lifting}
E^x\bigl(\lift{\phi}{\hilc}(T)\bigr)E_y=\phi(E^x T E_y)\qquad%
\forall\,T\in\mat{\hilc}{\opsp_1},\ x,y\in\hilc.
\end{equation}
The \emph{lifting} $\lift{\phi}{\hilc}$ is linear and $\hilc$ bounded,
with $\|\lift{\phi}{\hilc}\|\leq\|\phi\|_\hilc$ and
$\|\lift{\phi}{\hilc}\|_\cb\leq\|\phi\|_\cb$.
\end{thm}
\begin{proof}
If $\hilc$ is finite dimensional then
$\lift{\phi}{\hilc}=\phi\algten\id_{\bop{\hilc}{}}$; otherwise, let
$T\in\mat{\hilc}{\opsp_1}$ and apply Zorn's lemma to the collection of
pairs $(M,R)$ such that $M$ is a closed subspace of $\hilc$,
$R\in\mat{M}{\opsp_2}$ and $E^x R E_y=\phi(E^x T E_y)$ for all
$x$,~$y\in M$, ordered by setting $(M,R)\leq(N,S)$ whenever
$M\subseteq N$; for any such $(M,R)$, it is simple to verify that
\begin{equation}\label{eqn:cbineq}
\|R\|=\sup_L\|P_L R|_{\hilb_2\otimes L}\|=%
\sup_L\|(\phi\algten\id_{\bop{L}{}})(P_L T|_{\hilb_1\otimes L})\|\leq%
\|\phi\|_\cb\|T\|,
\end{equation}
where $P_L\in\bop{\hilb_i\otimes L}{}$ is the orthogonal projection
with range $\hilb_i\otimes L$ ($i=1$,~$2$) and the supremum is taken
over all finite-dimensional subspaces of $M$.

To see that $T\mapsto\lift{\phi}{\hilc}(T)$ has the desired
properties, note that uniqueness and linearity follow from
(\ref{eqn:lifting}). Proposition~\ref{prp:hsineq} and
(\ref{eqn:cbineq}) give the inequalities; if $n\geq1$ then
\[
\|(\lift{\phi}{\hilc})\algten\id_{\bop{\C^n}{}}\|=%
\|\lift{(\phi\algten\id_{\bop{\C^n}{}})}{\hilc}\|\leq%
\|\phi\algten\id_{\bop{\C^n}{}}\|_\cb\leq%
\|\phi\|_\cb.\qedhere
\]
\end{proof}

\begin{rem}\label{rmk:alglift}
Let $\phi\in\kbo{\hilc}{\opsp_1}{\opsp_2}$, where $\opsp_1$ and
$\opsp_2$ are operator spaces.
\begin{enumerate}
\item[(i)] The restriction 
$\lift{\phi}{\hilc}|_{\opsp_1\splten\bop{\hilc}{}}$ equals
$\phi\splten\id_{\bop{\hilc}{}}$.
\item[(ii)] If $\opsp_1$ and $\opsp_2$ are ultraweakly closed and
$\phi$ is ultraweakly continuous then
$\lift{\phi}{\hilc}=\phi\uwkten\id_{\bop{\hilc}{}}$. (Note that
$\lift{\phi}{\hilc}$ extends $\phi\splten\id_{\bop{\hilc}{}}$ even if
$\phi$ is not ultraweakly continuous. Neufang \cite[{\S}5]{Neu04} has
examined this phenomenon.)
\item[(iii)] The $\hilc$-bounded map
\[
\phi^\dagger:\opsp_1^\dagger\to\opsp_2^\dagger;\ a^*\mapsto\phi(a)^*
\]
is completely bounded if $\phi$ is and
$\lift{\phi^\dagger}{\hilc}=(\lift{\phi}{\hilc})^\dagger$.
\end{enumerate}
\end{rem}

\gap
The family of maps $(\phi^{(m)})_{m\geq0}$ defined in the following
theorem is the \emph{quantum random walk} with
\emph{generator} $\phi$.

\begin{thm}\label{thm:qrw}
If $\phi\in\hilc\bopp\bigl(\opsp;\mat{\hilc}{\opsp}\big)$ then there
exists a unique family of maps
$\phi^{(m)}:\opsp\to\mat{\hilc^{\otimes m}}{\opsp}$ such that
$\phi^{(0)}=\id_\opsp$ and
\begin{equation}\label{eqn:come}
E^x \phi^{(m)}(a) E_y=\phi^{(m-1)}(E^x \phi(a) E_y)%
\qquad\forall\,x,y\in\hilc,\ a\in\opsp,\ m\geq1.
\end{equation}
These maps are necessarily linear and $\hilc$ bounded, and are
completely bounded if $\phi$ is, with
$\|\phi^{(m)}\|_\hilc\leq\|\phi\|_\hilc^m$ and
$\|\phi^{(m)}\|_\cb\leq\|\phi\|_\cb^m$ for all $m\geq1$.
\end{thm}
\begin{proof}
Let $\phi^{(0)}:=\id_\opsp$ and
$\phi^{(m+1)}:=(\lift{\phi}{\hilc^{\otimes m}})\circ\phi^{(m)}$ for
all $m\in\Z_+$, with the spaces
$\bigl(\mat{\hilc}{\opsp}\bigr)\matten\bop{\hilc^{\otimes m}}{}$ and
$\mat{\hilc^{\otimes m+1}}{\opsp}$ identified in the natural manner;
uniqueness and $\hilc$-boundedness are clear. For the first
inequality, suppose $\dim\hilc<\infty$ and let $\{e_i:i\in I\}$ be an
orthonormal basis for $\hilc$. If $m\geq1$, $a\in\opsp$ and
$x$,~$y\in\hilb\otimes\hilc^{m+1}$ then, with $\bp=(p_1,\ldots,p_m)$
and $\bq=(q_1,\ldots,q_m)$,
\begin{align*}
|\langle x,\phi^{(m+1)}(a)y\rangle|&\leq%
\biggl|\sum_{\bp,\bq\in I^m}\langle E^{e_{p_1}\otimes\cdots%
\otimes e_{p_m}}x,\phi\circ\phi^{p_1}_{q_1}\circ\cdots%
\circ\phi^{p_m}_{q_m}(a)%
E^{e_{q_1}\otimes\cdots\otimes e_{q_m}}y\rangle\biggr|\\
&\leq\|\phi\|^{m+1}\|a\|%
\sum_{\bp\in I^m}\|E^{e_{p_1}\otimes\cdots\otimes e_{p_m}}x\|%
\sum_{\bq\in I^m}\|E^{e_{q_1}\otimes\cdots\otimes e_{q_m}}y\|\\
&\leq\|\phi\|^{m+1}\|a\|\,(\dim\hilc)^m\|x\|\,\|y\|,
\end{align*}
where $\phi^p_q:=E^{e_p}\phi(\cdot)E_{e_q}$ for all
$p$,~$q\in I$. Hence
$\|\phi^{(m+1)}\|\leq\|\phi\|^{m+1}(\dim\hilc)^m$, as required. The
second inequality is immediate.
\end{proof}

\begin{rem}\label{rmk:walkadj}
If $\phi\in\hilc\bopp\bigl(\opsp;\mat{\hilc}{\opsp}\bigr)$ and
$m\in\Z_+$ then $(\phi^\dagger)^{(m)}=(\phi^{(m)})^\dagger$.
\end{rem}

\begin{prop}\label{prp:hom}
Let $\csa\subseteq\bop{\hilb}{}$ be a $C^*$~algebra.
\begin{enumerate}
\item[(i)] If $\phi:\csa\to\csa\splten\bop{\hilc}{}$ is a
$*$-homomorphism then, for all $m\in\Z_+$, so is
the map $\phi^{(m)}:\csa\to\csa\splten\bop{\hilc}{}^{\splten m}$.
\item[(ii)] If $\csa$ is a von~Neumann algebra and
$\phi:\csa\to\csa\uwkten\bop{\hilc}{}$ is a normal (\ie ultraweakly
continuous) $*$-homomorphism then, for all $m\in\Z_+$, so is
$\phi^{(m)}:\csa\to\csa\uwkten\bop{\hilc^{\otimes m}}{}$.
\end{enumerate}
In both cases, the maps $(\phi^{(m)})_{m\geq0}$ are completely
bounded.
\end{prop}
\begin{proof}
As a $C^*$-algebra $*$-homomorphism is contractive
\cite[Theorem~2.1.7]{Mur90}, so completely contractive, induction and
Remark~\ref{rmk:alglift} give the result.
\end{proof}

\begin{exmp}{[\Cf\cite{AtP06,Bru02,Gou04}.]}
If $\agb$ is a von~Neumann algebra, the operator
$U\in\agb\uwkten\bop{\hilc}{}$ is unitary and the normal unital
$*$-homomorphism
\[
\phi:\agb\to\agb\uwkten\bop{\hilc}{};\ a\mapsto %
U^*(a\otimes\id_\hilc)U
\]
then the quantum random walk $(\phi^{(m)})_{m\geq0}$ admits the
following physical interpretation. Elements of the algebra $\agb$
describe the configuration of a system, which interacts periodically
with a series of identical particles whose configurations are
described by elements of $\bop{\hilc}{}$; the interaction between the
system and an individual particle is given by the map $\phi$. If the
system is initially in configuration $a$, the operator $\phi^{(m)}(a)$
represents the combined configuration of the system and the first $m$
particles with which it has interacted, with the particles `moving to
the right' in the product
$\bop{\hilc}{}\uwkten\cdots\uwkten\bop{\hilc}{}=%
\bop{\hilc^{\otimes m}}{}$.
\end{exmp}

\begin{prop}\label{prp:liftcng}
If $\phi_n\in\kbo{\hilc}{\opsp_1}{\opsp_2}$ is such that
$\lift{\phi_n}{\hilc}\to0$ strongly as $n\to\infty$, \ie
\[
\lim_{n\to\infty}\|\lift{\phi_n}{\hilc}(b)\|=0\qquad%
\forall\,b\in\mat{\hilc}{\opsp_1}
\]
then $\lift{\phi_n}{\hilc^{\otimes m}}\to0$ strongly as $n\to\infty$
for all $m\in\Z_+$.
\end{prop}
\begin{proof}
For $m=0$, choose a unit vector $e\in\hilc$ and note that, as
$n\to\infty$,
\[
\|\phi_n(a)\|=%
\|\lift{\phi_n}{\hilc}(a\otimes\dyad{e}{e})\|\to0%
\qquad\forall\,a\in\opsp_1,
\]
where $\dyad{e}{e}\in\bop{\hilc}{}$ is the orthogonal projection onto
$\C e$.

If $\hilc$ is finite dimensional then the result now follows, since
$\mat{\hilc^{\otimes m}}{\opsp_1}$ is the linear span of simple
tensors. Otherwise, note that $\hilc^{\otimes m}$ is isomorphic to
$\hilc$ for all $m\geq1$ and if $U:\hilc_1\to\hilc_2$ is a unitary
operator then
\begin{equation}\label{eqn:isom}
(\id_\hilb\otimes U)%
\bigl(\lift{\phi_n}{\hilc_1}(\cdot)\bigr)(\id_\hilb\otimes U)^*=%
\lift{\phi_n}{\hilc_2}\bigl((\id_\hilb\otimes U)\cdot%
(\id_\hilb\otimes U)^*\bigr),
\end{equation}
assuming the space $\opsp_1\subseteq\bop{\hilb}{}$.
\end{proof}

\begin{lem}\label{lem:approxgen}
Suppose
$\phi_n$,~$\phi\in\hilc\bopp\bigl(\opsp;\mat{\hilc}{\opsp}\bigr)$ and
let $m\in\Z_+$. If $\lift{\phi_n}{\hilc}\to\lift{\phi}{\hilc}$
strongly as $n\to\infty$ then $\phi_n^{(m)}\to\phi^{(m)}$ strongly; if
$\phi_n\to\phi$ in $\hilc$ norm or $\cb$ norm, with $\phi_n$ and
$\phi$ completely bounded in the latter case, then
$\phi_n^{(m)}\to\phi^{(m)}$ in the same sense.
\end{lem}
\begin{proof}
As $\phi^{(m+1)}:=\lift{\phi}{\hilc^{\otimes m}}\circ\phi^{(m)}$ for
all $m\in\Z_+$, the first claim follows by induction from
Proposition~\ref{prp:liftcng}, together with the estimate
\begin{multline*}
\|(\phi_n^{(m+1)}-\phi^{(m+1)})(a)\|\leq%
\|\lift{\phi_n}{\hilc^{\otimes m}}\|\,%
\|(\phi_n^{(m)}-\phi^{(m)})(a)\|\\
+\|\bigl(\lift{(\phi_n-\phi)}{\hilc^{\otimes m}}\bigr)\circ%
\phi^{(m)}(a)\|
\end{multline*}
and the principle of uniform boundedness. As
$\|\cdot\|_{\hilc^{\otimes p}}\leq\|\cdot\|_{\hilc^{\otimes q}}$ if
$p\leq q$ and
\[
\|\phi_n^{(m+1)}-\phi^{(m+1)}\|_\hilc\leq%
\|\phi_n\|_{\hilc^{\otimes m}}\|\phi_n^{(m)}-\phi^{(m)}\|_\hilc+%
\|\phi_n-\phi\|_{\hilc^{\otimes m}}\|\phi^{(m)}\|_\hilc,
\]
induction and Theorem~\ref{thm:qrw} yield the inequality
\[
\|\phi_n^{(m+1)}-\phi^{(m+1)}\|_\hilc\leq%
\|\phi_n-\phi\|_{\hilc^{\otimes m}}%
\bigl(\|\phi_n\|_{\hilc^{\otimes m}}+\|\phi\|_{\hilc^{\otimes m}}\bigr)^{m}.
\]
Since $\|\cdot\|_\hilc$ and $\|\cdot\|_{\hilc^{\otimes m}}$ are
equivalent, the result follows; the same working applies to the
completely bounded case.
\end{proof}

\gap
The following result was brought to our attention by
Skalski~\cite{Ska08}.

\begin{lem}\label{lem:cbequiv}
Let $\hilc$ be infinite dimensional and let
$\phi_n\in\cbo{\opsp_1}{\opsp_2}$. Then $\lift{\phi_n}{\hilc}\to0$
strongly as $n\to\infty$ if and only if $\phi_n\to0$ in $\cb$ norm.
\end{lem}
\begin{proof}
Since $\|\lift{\phi_n}{\hilc}\|\leq\|\phi_n\|_\hilc$, one direction is
immediate. For the other, suppose $\lift{\phi_n}{\hilc}\to0$ strongly
as $n\to\infty$, so that $\lift{\phi_n}{\hilc_1}\to0$ strongly if
$\hilc_1$ has the same dimension as $\hilc$, by (\ref{eqn:isom}), and
$\|\phi_n\|_\cb\not\to0$. There exists $\eps>0$ such that, passing to
a subsequence if necessary, $\|\phi_n\|_\cb>\eps$ for all $n\geq1$:
choose $m_n\geq1$ and $c_n\in\opsp_1\algten\bop{\C^{m_n}}{}$ such that
$\|c_n\|=1$ and $\|\phi_n\algten\id_{\bop{\C^{m_n}}{}}(c_n)\|>\eps$.
Let $\hilc_1:=\hilc\oplus\bigoplus_{n=1}^\infty \C^{m_n}$ and recall
that if $T_{(n)}\in\bop{\hilb_{(n)}}{}$ for all $n\geq0$ then
\[
\oplus_{n=0}^\infty T_{(n)}:\bigoplus_{n=0}^\infty \hilb_{(n)}%
\to\bigoplus_{n=0}^\infty \hilb_{(n)};\ %
(x_n)_{n=0}^\infty\mapsto(T_{(n)}x_n)_{n=0}^\infty
\]
has norm equal to $\sup_{n\geq0}\|T_{(n)}\|$. Hence
$c:=\oplus_{n=0}^\infty c_n$ has norm $1$, where $c_0:=0$, and it is
readily verified that
$c\in\opsp_1\matten\bop{\hilc_1}{}$. Furthermore,
$\phi_p\matten\id_{\bop{\hilc_1}{}}(c)=%
\oplus_{n=1}^\infty\phi_p\algten\id_{\bop{\C^{m_n}}{}}(c_n)$, so
$\|\phi_p\matten\id_{\bop{\hilc_1}{}}(c)\|>\eps$ for all
$p\geq1$. This is the desired contradiction.
\end{proof}

\section{Toy Fock space and quantum stochastic integrals}\label{sec:toy}

\begin{notation}
Let $\mul$ be a complex Hilbert space (the \emph{multiplicity space})
and let $\mmul:=\C\oplus\mul$ be its one-dimensional
extension. Elements of $\mmul$ will be thought of as column vectors,
with the first entry a complex number and the second a vector in
$\mul$; if $x\in\mul$ then $\widehat{x}:=\mvec{1}{x}$. This
decomposition will be used to write various operators as two-by-two
matrices.
\end{notation}

\begin{defn}
\emph{Toy Fock space} is the countable tensor product
\[
\bebe:=\bigotimes_{n=0}^\infty\mmul_{(n)}
\]
with respect to the stabilising sequence
$\bigl(\vac_{(n)}:=\mvec{1}{0}\bigr)_{n\geq0}$, where
$\mmul_{(n)}:=\mmul$ for each $n$; the subscript $(n)$ is used here
and below to indicate the relevant copy. (For information on infinite
tensor products of Hilbert spaces, see, for example
\cite[Exercise~11.5.29]{KaR97}.) For all $n\in\Z_+$, let
\[
\bebe_{n)}:=\bigotimes_{m=0}^{n-1}\mmul_{(m)}\quad\mbox{and}\quad%
\bebe_{[n}:=\bigotimes_{m=n}^\infty\mmul_{(m)},
\]
where $\bebe_{0)}:=\C$, so that
$\bebe=\bebe_{n)}\otimes\bebe_{[n}$; this is the analogue of the
continuous tensor-product structure of Boson Fock space.
\end{defn}

\begin{notation}
Let $\fock$ denote Boson Fock space over $\elltwo$, the complex
Hilbert space of $\mul$-valued, square-integrable functions on the
half line, and let $\evecs:=\lin\{\evec{f}:f\in\elltwo\}$ be the
subspace of $\fock$ spanned by exponential vectors, where $\evec{f}$
denotes the exponential vector corresponding to $f$. Let $\ini$ be a
complex Hilbert space (the \emph{initial space}) and let
$\bbebe:=\ini\otimes\bebe$, $\ffock:=\ini\otimes\fock$,
$\eevecs:=\ini\algten\evecs$ \etc.
\end{notation}

\begin{thm}\label{thm:embed}
For all $h>0$ there exists a unique co-isometry $D_h:\ffock\to\bbebe$
such that
\[
D_h u\evec{f}=u\otimes\bigotimes_{n=0}^\infty\widehat{f(n;h)}%
\qquad\forall\,u\in\ini,\ f\in\elltwo,
\]
where the tensor sign between components of simple tensors in $\ffock$
is omitted and
\[
f(n;h):=h^{-1/2}\int_{n h}^{(n+1)h}f(t)\std t\qquad\forall\,n\in\Z_+.
\]
Moreover, $D_h^*D_h\to\id_\ffock$ in the strong operator topology as
$h\to0$.
\end{thm}
\begin{proof}
See \cite[Definition~2.2, Theorem~2.1 and Notation~4.1]{Blt07}.
\end{proof}

\begin{defn}
Given a Hilbert space $\hilb$, an \emph{$\hilb$~process}
$X=(X_t)_{t\in\R_+}$ is a weakly measurable family of linear
operators with common domain $\hilb\algten\evecs$, \ie
\[
X_t\in\lop{\hilb\algten\evecs}{\hilb\otimes\fock}\qquad\forall\,t\in\R_+
\]
and $t\mapsto\langle u\evec{f},X_t v\evec{g}\rangle$ is measurable for
all $u$,~$v\in\hilb$ and $f$,~$g\in\elltwo$.

An $\hilb$~process $X$ is \emph{vacuum adapted} if
\[
\langle u\evec{f},X_t v\evec{g}\rangle=%
\langle u\evec{\indf{[0,t[}f},X_t v\evec{\indf{[0,t[}g}\rangle
\]
for all $t\in\R_+$, $u$,~$v\in\hilb$ and
$f$,~$g\in\elltwo$. Equivalently, the identity
$(\id_\hilb\otimes\expn_t)X_t(\id_\hilb\otimes\expn_t)=X_t$ holds for
all $t\in\R_+$, where $\expn_t\in\bop{\fock}{}$ is the second
quantisation of the multiplication operator $f\mapsto\indf{[0,t[}f$ on
$\elltwo$.

An $\hilb$~process $X$ is \emph{semi-vacuum-adapted} if
$(\id_\hilb\otimes\expn_t)X_t=X_t$ for all $t\in\R_+$; clearly every
vacuum-adapted process is semi-vacuum-adapted. (The modified integral
of Theorem~\ref{thm:vacsubo} below preserves semi-vacuum-adaptedness
but not vacuum adaptedness.)
\end{defn}

\begin{notation}
For all $m\geq1$ and $t\in\R_+$, let
\[
\Delta_m(t):=\{\bt:=(t_1,\ldots,t_m)\in{[0,t[}^m:t_1<\cdots<t_m\}%
\subseteq\R_+^m
\]
and, given $M\in\bop{\elltwo}{}$, define
$\nnabla^{M,m}\in\lopp\bigl(\eevecs;%
L^2(\Delta_m(t);\ini\otimes\mmul^{\otimes m}\otimes\fock)\bigr)$
such that
\[
\nnabla^{M,m}_\bt u\evec{f}:=%
\bigl(\nnabla^{M,m}u\evec{f}\bigr)(\bt):=%
u\otimes\widehat{M f}^{\otimes m}(\bt)\otimes\evec{f}
\]
for all $u\in\ini$ and $f\in\elltwo$, where
$\widehat{g}^{\otimes m}(\bt):=%
\widehat{g(t_1)}\otimes\cdots\otimes\widehat{g(t_m)}$ for all
$g\in\elltwo$ and $\bt\in\R_+^m$. For brevity, let
$\nnabla^m_\bt:=\nnabla^{M,m}_\bt$ when $M=\id_\elltwo$.
\end{notation}

\begin{thm}\label{thm:vaciter}
If $m\geq1$, $X\in\bop{\bbebe_{m)}}{}$ and $Y$ is a locally uniformly
bounded, vacuum-adapted $\C$~process there exists a unique
vacuum-adapted $\ini$~process $\vint^m(X\otimes Y)$, the
\emph{$m$-fold QS integral}, such that, for all $t\in\R_+$,
\begin{equation}\label{eqn:vinm}
\|\vint^m(X\otimes Y)_t\eta\|^2\leq c_t^{2m}%
\int_{\Delta_m(t)}\|(X\otimes Y_{t_1})%
\nnabla^m_\bt\eta\|^2\std\bt
\end{equation}
for all $\eta\in\eevecs$, where $c_t:=\sqrt{2\max\{t,1\}}$, and
\begin{multline}\label{eqn:viip}
\langle u\evec{f},\vint^m(X\otimes Y)_t v\evec{g}\rangle\\
=\int_{\Delta_m(t)}\langle u\otimes\widehat{f}^{\otimes m}(\bt),%
X\bigl(v\otimes\widehat{g}^{\otimes m}(\bt)\bigr)\rangle%
\langle\evec{f},Y_{t_1}\evec{g}\rangle\std\bt
\end{multline}
for all $u$,~$v\in\ini$ and $f$,~$g\in\elltwo$.
\end{thm}
\begin{proof}
See \cite[Theorem~3.2]{Blt07}.
\end{proof}

\begin{defn}
If $A$ is an ordered set and $m\in\Z_+$ then $A^{m,\uparrow}$ denotes
the collection of strictly increasing $m$-tuples of elements of $A$,
with $A^{0,\uparrow}:=\{\emptyset\}$. Given $h>0$ and $m\geq1$, let
\[
{[\bp h,(\bp+1)h[}:=%
\{(t_1,\ldots,t_m)\in\R_+^m:p_i h\leq t_i<(p_i+1)h\ %
\quad(i=1,\ldots,m)\}
\]
for all $\bp=(p_1,\ldots,p_m)\in\Z_+^{m,\uparrow}$ and, for all
$t\in\R_+$, let
\[
\Delta_m^h(t):=\bigcup_{\bp\in\{0,\ldots,n-1\}^{m,\uparrow}}%
{[\bp h,(\bp+1)h[}\qquad\mbox{if }t\in{[n h,(n+1)h[}.
\]
\end{defn}

\begin{thm}\label{thm:vacsubo}
Let $M\in\bop{\elltwo}{}$. If $m\geq1$, $h>0$,
$X\in\bop{\bbebe_{m)}}{}$ and $Y$ is a locally uniformly bounded,
semi-vacuum-adapted $\C$~process then there exists a unique
semi-vacuum-adapted $\ini$~process $\vint^m(X\otimes Y;M,h)$, the
\emph{modified $m$-fold QS integral of step size $h$}, such that,
for all $t\in\R_+$,
\begin{equation}\label{eqn:subonorm}
\|\vint^m(X\otimes Y;M,h)_t\eta\|^2\leq c_t^{2m}%
\int_{\Delta_m^h(t)}\|(X\otimes Y_{t_1})%
\nnabla^{M,m}_\bt\eta\|^2\std\bt
\end{equation}
for all $\eta\in\eevecs$ and
\begin{multline*}
\langle u\evec{f},\vint^m(X\otimes Y;M,h)_t v\evec{g}\rangle\\
=\int_{\Delta_m^h(t)}\langle u\otimes\widehat{f}^{\otimes m}%
(\bt),X\bigl(v\otimes\widehat{M g}^{\otimes m}(\bt)\bigr)\rangle%
\langle\evec{f},Y_{t_1}\evec{g}\rangle\std\bt
\end{multline*}
for all $u$,~$v\in\ini$ and $f$,~$g\in\elltwo$.
\end{thm}
\begin{proof}
See \cite[Theorem~3.3]{Blt07}.
\end{proof}

\begin{defn}
For all $h>0$, let
$\Xi_h:=\left[\begin{smallmatrix}
h^{-1/2}&0\\0&\id_\mul
\end{smallmatrix}\right]\in\bop{\mmul}{}$
and, for all $m\in\Z_+$, define
\[
\scale_h(X):=(\id_\ini\otimes\Xi_h^{\otimes m})X%
(\id_\ini\otimes\Xi_h^{\otimes m})\qquad%
\mbox{if }X\in\bop{\bbebe_{m)}}{},
\]
so $\scale_h(X):=X$ if $X\in\bop{\ini}{}$. The mapping $\scale_h$
is $*$-linear but not, in general, multiplicative.
\end{defn}

\begin{prop}\label{prp:scale}
If $h>0$ and $\phi\in\mmul\bopp\bigl(\opsp;\mat{\mmul}{\opsp}\bigr)$
then the scaled map
$\scale_h(\phi):a\mapsto\scale_h\bigl(\phi(a)\bigr)$ is $\mmul$
bounded from $\opsp$ to $\mat{\mmul}{\opsp}$, and is completely
bounded if $\phi$ is. Furthermore,
\begin{equation}\label{eqn:scid}
\scale_h(\phi)^{(n)}:=\bigl(\scale_h(\phi)\bigr)^{(n)}=%
\scale_h(\phi^{(n)})\qquad\forall\,n\in\Z_+.
\end{equation}
\end{prop}
\begin{proof}
Note first that
$E^x\scale_h\bigl(\phi(a)\bigr)E_y=%
E^{\Xi_h x}\phi(a)E_{\Xi_h y}\in\opsp$ for all
$x$,~$y\in\mmul$, so
$\scale_h\bigl(\phi(a)\bigr)\in\mat{\mmul}{\opsp}$. Next,
\begin{align*}
\|(\scale_h(\phi)\algten\id_{\bop{\C^n}{}})(T)\|&=%
\|(\id_\ini\otimes\Xi_h\otimes\id_{\C^n})%
(\phi\algten\id_{\bop{\C^n}{}})(T)%
(\id_\ini\otimes\Xi_h\otimes\id_{\C^n})\|\\
&\leq\max\{1,h^{-1}\}\|\phi\algten\id_{\bop{\C^n}{}}\|\,\|T\|
\end{align*}
for all $T\in\opsp\algten\bop{\C^n}{}$, hence
$\|\scale_h(\phi)\|_\mmul\leq\max\{1,h^{-1}\}\|\phi\|_\mmul$ \etc. As
(\ref{eqn:scid}) holds for $n=0$ and $n=1$, the general case follows
by induction.
\end{proof}

\begin{defn}\label{def:vhm}
For all $m\in\Z_+$ and $\bp\in\Z_+^{m,\uparrow}$, let
$\vhm_\bp:\bop{\bbebe_{m)}}{}\to\bop{\bbebe}{}$ be the normal
$*$-homomorphism such that
\[
X\otimes B_1\otimes\cdots\otimes B_m\mapsto X\otimes%
\Delta^\perp_{[0,p_1)}\otimes B_1\otimes%
\Delta^\perp_{[p_1+1,p_2)}\otimes\cdots\otimes B_m%
\otimes\Delta^\perp_{[p_m+1},
\]
on the right-hand side of which $B_n$ acts on $\mmul_{(p_n)}$ for
$n=1,\ldots,m$ and the vacuum projection
$\Delta^\perp:=\dyad{\vac}{\vac}:x\mapsto\langle\vac,x\rangle\vac$
acts on $\mmul_{(q)}$ whenever $q\not\in\{p_1,\ldots,p_m\}$. (In
particular, $\vhm_\emptyset(X)=X\otimes\Delta^\perp_{[0}$ for all
$X\in\bop{\ini}{}$.)
\end{defn}

\begin{notation}
If $h>0$ then
\[
P_{(h)}f:=\sum_{n=0}^\infty\frac{1}{h}\int_{n h}^{(n+1)h}f(t)\std t%
\ \indf{[n h,(n+1)h[}\qquad\forall\,f\in\elltwo,
\]
so $P_{(h)}\in\bopp\bigl(\elltwo\bigr)$ is an orthogonal projection and
$P_{(h)}\to\id_\elltwo$ strongly as $h\to0$.
\end{notation}

\gap
The following processes, which are fundamental for the approximation
results in Section~\ref{sec:wlk}, are a variation on the toy integrals
of \cite{Blt07}. They appear naturally when a quantum random walk is
re-written using a Wiener--It\^o decomposition
(Proposition~\ref{prp:qrwwi}).

\begin{prop}\label{prp:vdi}
If $m\geq1$, $h>0$, $X\in\bop{\bbebe_{m)}}{}$ and $t\in\R_+$ then
\[
\vint^m(X\otimes\expn_\Vac;P_{(h)},h)_t=%
\sum_{\bp\in\Z_+^{m,\uparrow}}\tfn{p_m+1\leq t/h}%
D_h^*\vhm_\bp\bigl(\scale_{h^{-1}}(X)\bigr)D_h,
\]
where $(\expn_\Vac)_s:=\expn_0$ for all $s\in\R_+$.
\end{prop}
\begin{proof}
This is a straightforward exercise (\cf\cite[Theorem~5.1]{Blt07}).
\end{proof}

\section{Vacuum cocycles}\label{sec:coc}

\begin{defn}
Let $\opsp\subseteq\bop{\ini}{}$ be an operator space and
$\psi:\opsp\to\mat{\mmul}{\opsp}$ a linear map. If there exists a
family of operators
\[
j_t\in\lopp\bigl(\opsp;\lop{\eevecs}{\ffock}\bigr)\qquad(t\in\R_+)
\]
such that $\bigl(j_t(a)\bigr)_{t\in\R_+}$ is a vacuum-adapted
$\ini$~process for all $a\in\opsp$, with
\begin{equation}\label{eqn:flow}
\langle u\evec{f},j_t(a) v\evec{g}\rangle=\langle u,a v\rangle+%
\int_0^t\langle u\evec{f},j_s(E^{\widehat{f(s)}}%
\psi(a)E_{\widehat{g(s)}}) v\evec{g}\rangle\std s
\end{equation}
for all $u$,~$v\in\ini$, $f$,~$g\in\elltwo$ and $t\in\R_+$, then $j$
is a type of mapping process called a \emph{vacuum flow} with
\emph{generator} $\psi$. The flow is \emph{weakly regular} if
\[
a\mapsto E^{\evec{f}}j_t(a) E_{\evec{g}}\in%
\bopp\bigl(\opsp;\bop{\ini}{}\bigr)%
\qquad\forall\,f,g\in\elltwo
\]
with norm locally uniformly bounded as a function of $t$. (It is part
of the definition of weak regularity that
$E^{\evec{f}}j_t(a)E_{\evec{g}}\in\bop{\ini}{}$ for all $a\in\opsp$;
this need not hold \textit{a priori}.)
\end{defn}

\begin{prop}
There is at most one weakly regular vacuum flow with generator $\psi$
for any given $\psi\in\bopp\bigl(\opsp;\mat{\mmul}{\opsp}\bigr)$.
\end{prop}
\begin{proof}
If $j$ and $j'$ are two such then fix $u$,~$v\in\ini$,
$f$,~$g\in\elltwo$, $t\in\R_+$ and $a\in\opsp$; for all $n\geq1$,
\begin{align*}
R&:=\langle u\evec{f},(j-j')_t(a) v\evec{g}\rangle\\
&\phantom{:}=\int_{\Delta_n(t)}\langle u,E^{\evec{f}}(j-j')_{t_1}%
\bigl(\psi^{\widehat{f(t_1)}}_{\widehat{g(t_1)}}\circ\cdots%
\circ\psi^{\widehat{f(t_n)}}_{\widehat{g(t_n)}}(a)\bigr)%
E_{\evec{g}}v\rangle\std\bt,
\end{align*}
where $\psi^x_y:=E^x\psi(\cdot)E_y$ for all $x$,~$y\in\mmul$. If
if $C_t$ denotes the norm of $E^{\evec{f}}j_t(\cdot)E_{\evec{g}}$ and
similarly for $C'_t$, then
\[
|R|\leq\|u\|\sup_{0\leq s\leq t}\{C_s+C'_s\}\,\|\psi\|^n\|a\|%
\Bigl(\int_0^t\bigl\|\widehat{f(s)}\bigr\|\,%
\bigl\|\widehat{g(s)}\bigr\|\std s\Bigr)^n\|v\|/n!,
\]
which tends to $0$ as $n\to\infty$.
\end{proof}

\gap
The Wiener--It\^o decomposition given by the following theorem is at
the heart of Theorem~\ref{thm:main}, the main result of
Section~\ref{sec:wlk}. There is an analogous version for any quantum
random walk (Proposition~\ref{prp:qrwwi}) and each iterated integral
appearing below is the limit of a corresponding discrete version.

\begin{thm}\label{thm:regcoc}
If $\psi\in\mmul\bopp\bigl(\opsp;\mat{\mmul}{\opsp}\bigr)$ and
$(\psi^{(m)})_{m\geq0}$ is the quantum random walk of
Theorem~\ref{thm:qrw} then
\[
j^\psi_t(a):=a\otimes\expn_0+%
\sum_{m=1}^\infty\vint^m(\psi^{(m)}(a)\otimes\expn_\Vac)_t
\]
is strongly convergent on $\eevecs$ for all $a\in\opsp$ and
$t\in\R_+$; the process $j^\psi$ is the weakly regular vacuum flow
with generator $\psi$.
\end{thm}
\begin{proof}{[\Cf\cite[Theorem~4.3]{Lin05}.]}
By Theorem~\ref{thm:qrw}, $\|\psi^{(m)}\|\leq\|\psi\|^m_\mmul$, so
(\ref{eqn:vinm}) yields the inequality
\begin{equation}\label{eqn:qinm}
\|\vint^m\bigl(\psi^{(m)}(a)\otimes\expn_\Vac\bigr)_t %
u\evec{f}\|\leq c_t^m\|\psi\|^m_\mmul\|a\|\,%
\|\indf{[0,t[}\widehat{f}\|^m\|u\|/\sqrt{m!},
\end{equation}
which proves convergence. Next, (\ref{eqn:viip}) and (\ref{eqn:come})
imply that
\begin{multline*}
\langle u\evec{f},%
\vint^m(\psi^{(m)}(a)\otimes\expn_\Vac)_t v\evec{g}\rangle\\
=\int_0^t\langle u\evec{f},%
\vint^{m-1}(\psi^{(m-1)}(E^{\widehat{f(s)}}%
\psi(a)E_{\widehat{g(s)}})\otimes\expn_\Vac)_s %
v\evec{g}\rangle\std s,
\end{multline*}
where $\vint^0(X\otimes Y)_s:=X\otimes Y_s$. This identity, with
induction, gives (\ref{eqn:flow}) and vacuum adaptedness; regularity
follows from the estimate (\ref{eqn:qinm}).
\end{proof}

\begin{rem}\label{rmk:cbcoc}
Let $\psi\in\mmul\bopp\bigl(\opsp;\mat{\mmul}{\opsp}\bigr)$, $m\geq1$,
$t\in\R_+$ and $f\in\elltwo$; the inequality (\ref{eqn:qinm}) implies
that
\[
\vint^m(\psi^{(m)}(\cdot)\otimes\expn_\Vac)_t E_{\evec{f}}%
\in\bopp\big(\opsp;\bop{\ini}{\ffock}\bigr).
\]
If $n\geq1$ is fixed and the unitary operators
\begin{alignat*}{2}
U_1:\ &\ini\otimes\fock\otimes\C^n\ &\to&\ \ini\otimes\C^n\otimes\fock\\
\mbox{and}\quad U_{2,p}:\ &\ini\otimes\mmul^{\otimes p}\otimes\C^n\ %
&\to&\ \ini\otimes\C^n\otimes\mmul^{\otimes p}
\end{alignat*}
act by exchanging the last two tensor components then
\begin{equation}\label{eqn:vicb}
U_1(\vint^m(\psi^{(m)}(\cdot)\otimes\expn_\Vac)_t%
\algten\id_{\bop{\C^n}{}})=%
\vint^m(\Psi^{(m)}(\cdot)\otimes\expn_\Vac)_t E_{\evec{f}},
\end{equation}
where
$U_{2,1}^*\Psi(\cdot)U_{2,1}=\psi\algten\id_{\bop{\C^n}{}}$. (This
holds because
\[
U_{2,m}^*\Psi^{(m)}(\cdot)U_{2,m}=\psi^{(m)}\algten\id_{\bop{\C^n}{}},
\]
which may be proved by induction on $m$.) In particular, if $\psi$ is
completely bounded then so is
$\vint^m(\psi^{(m)}(\cdot)\otimes\expn_\Vac)_t E_{\evec{f}}$; as
$\|\Psi\|_\cb\leq\|\psi\|_\cb$, it follows that
\begin{equation}\label{eqn:vicbnm}
\|\vint^m(\psi^{(m)}(\cdot)\otimes\expn_\Vac)_t %
E_{\evec{f}}\|_\cb\leq %
c_t^m\|\indf{[0,t[}\widehat{f}\|^m\|\psi\|^m_\cb/\sqrt{m!}.
\end{equation}
This (and the previous theorem) implies that
$j^\psi_t(\cdot)E_{\evec{f}}\in%
\mmul\bopp\bigl(\opsp;\bop{\ini}{\ffock}\bigr)$
and is completely bounded if $\psi$ is.
\end{rem}

\begin{prop}\label{prp:feller}
If $\psi\in\mmul\bopp\bigl(\opsp;\mat{\mmul}{\opsp}\bigr)$ then
$j^\psi$ is a Feller cocycle on $\opsp$ in the sense of
\emph{\cite{LiW00}}: for all $a\in\opsp$, $s$,~$t\in\R_+$ and
$f$,~$g\in\elltwo$,
\begin{align}
&E^{\evec{0}}j^\psi_0(a) E_{\evec{0}}=a,\nonumber\\
&E^{\evec{\indf{[0,t[}f}}j^\psi_t(a)%
E_{\evec{\indf{[0,t[}g}}\in\opsp\quad\mbox{and}\label{eqn:ccc2}\\
&E^{\evec{\indf{[0,s+t[}f}}j^\psi_{s+t}(a)%
E_{\evec{\indf{[s+t[}g}}\nonumber\\
&\quad=%
E^{\evec{\indf{[0,s[}f}}j^\psi_s(E^{\evec{\indf{[0,t[}f(\cdot+s)}}%
j^\psi_t(a)E_{\evec{\indf{[0,t[}g(\cdot+s)}})%
E_{\evec{\indf{[0,s[}g}}.\label{eqn:ccc3}
\end{align}
\end{prop}
\begin{proof}{[\Cf\cite[Theorem~5.1]{LiW00}.]}
Since $t\mapsto j^\psi_t(a)$ is vacuum adapted for all $a\in\opsp$,
\[
E^{\evec{\indf{[0,t[}f}}j^\psi_t(a)E_{\evec{\indf{[0,t[}g}}=%
E^{\evec{f}}j^\psi_t(a)E_{\evec{g}}%
\qquad\forall\,f,g\in\elltwo,\ t\in\R_+;
\]
if $j^{f,g}_t\in\bopp\bigl(\opsp;\bop{\ini}{}\bigr)$ is the mapping
$a\mapsto E^{\evec{f}}j^\psi_t(a)E_{\evec{g}}$ then it is required to
show that
\[
j^{f,g}_0(a)=a,\qquad j^{f,g}_t(a)\in\opsp\qquad\mbox{and}\qquad %
j^{f,g}_{s+t}=j^{f,g}_s\circ j^{f(\cdot+s),g(\cdot+s)}_t
\]
for all $f$,~$g\in\elltwo$, $a\in\opsp$ and $s$,~$t\in\R_+$. The first
of these is immediately verified. Next, let $f$ and $g$ be
$\mul$-valued step functions on $\R_+$ subordinate to the partition
$\{\tau_0<\tau_1<\cdots\}$, with values $f_n$ and $g_n$ respectively
on ${[\tau_n,\tau_{n+1}[}$ for all $n\in\Z_+$. If $u$,~$v\in\ini$ and
$t\in[0,\tau_1[$ then (\ref{eqn:flow}) implies that
\[
\langle u,j^{f,g}_t(a) v\rangle=\langle u,a v\rangle+%
\int_0^t\langle u,%
j^{f,g}_s\bigl(\psi^{\widehat{f_0}}_{\widehat{g_0}}(a)\bigr)v\rangle\std s,
\]
where $\psi^x_y\in\bop{\opsp}{}$ is defined by setting
$\psi^x_y(a):=E^x\psi(a)E_y$ for all
$x$,~$y\in\smash[t]{\mmul}$. If~$X$ and $Y$ are Banach spaces,
$A\in\bop{X}{}$ and $F:{[0,t_0[}\to\bopp\bigl(X;\bop{Y}{}\bigr)$ is
such that $t\mapsto F_t(x)y$ is bounded, measurable and satisfies the
integral equation
\[
F_t(x)y=F_0(x)y+\int_0^t\bigl((F_s\circ A)(x)\bigr)y\std s%
\qquad\forall\,t\in{[0,t_0[},\ x\in X,\ y\in Y,
\]
all in the weak sense, then $F_t=F_0\circ\exp(t A)$ on
${[0,t_0[}$. Hence
$\smash[t]{j^{f,g}_t=\exp(t\psi^{\widehat{f_0}}_{\widehat{g_0}})}$ for all
$t\in{[0,\tau_1[}$ (and for $t=\tau_1$, by continuity); in particular,
$j^{f,g}_t(a)\in\opsp$. More generally, if $n\geq1$ then
\[
\langle u,j^{f,g}_t(a)v\rangle=%
\langle u,j^{f,g}_{\tau_n}(a)v\rangle+\int_{\tau_n}^t\langle u,%
j^{f,g}_s\bigl(\psi^{\widehat{f_n}}_{\widehat{g_n}}(a)\bigr)v%
\rangle\std s
\]
for all $t\in{[\tau_n,\tau_{n+1}[}$ and so, for such $t$, we have the
semigroup decomposition
\[
j^{f,g}_t=j^{f,g}_{\tau_n}\circ%
\exp\bigl((t-\tau_n)\psi^{\widehat{f_n}}_{\widehat{g_n}}\bigr)=%
\exp\bigl((\tau_1-\tau_0)%
\psi^{\widehat{f_0}}_{\widehat{g_0}}\bigr)\circ\cdots\circ%
\exp\bigl((t-\tau_n)\psi^{\widehat{f_n}}_{\widehat{g_n}}\bigr).
\]
If $s\in{[\tau_m,\tau_{m+1}[}$ and $s+t\in{[\tau_n,\tau_{n+1}[}$ then
(assuming without loss of generality that $n>m$)
\begin{align*}
j^{f,g}_{s+t}&=\exp\bigl((\tau_1-\tau_0)%
\psi^{\widehat{f_0}}_{\widehat{g_0}}\bigr)\circ\cdots\circ%
\exp\bigl((s-\tau_m)\psi^{\widehat{f_m}}_{\widehat{g_m}}\bigr)\\
&\quad\circ\exp\bigl((\tau_{m+1}-s)%
\psi^{\widehat{f_m}}_{\widehat{g_m}}\bigr)\circ\cdots\circ%
\exp\bigl((t-(\tau_n-s))\psi^{\widehat{f_n}}_{\widehat{g_n}}\bigr)\\
&=j^{f,g}_s\circ j^{f(\cdot+s),g(\cdot+s)}_t,
\end{align*}
where the decomposition of $j^{f(\cdot+s),g(\cdot+s)}$ is taken with
respect to the partition $\{0<\tau_{m+1}-s<\tau_{m+2}-s<\cdots\}$. The
observation that $\mul$-valued step functions are dense in $\elltwo$
completes the proof.
\end{proof}

\begin{rem}
The flow $j^\psi$ in Proposition~\ref{prp:feller} is vacuum adapted,
so the identities (\ref{eqn:ccc2}) and (\ref{eqn:ccc3}) take the
simpler form
\begin{align*}
E^{\evec{f}}j^\psi_t(a)E_{\evec{g}}&\in\opsp\\
\mbox{and}\quad E^{\evec{f}}j_{s+t}(a)E_{\evec{g}}&=%
E^{\evec{f}}j^\psi_s(E^{\evec{f(\cdot+s)}}j_t(a)%
E_{\evec{g(\cdot+s)}})E_{\evec{g}}.
\end{align*}
\end{rem}

\begin{rem}\label{rmk:invol2}
If $\psi\in\mmul\bopp\bigl(\opsp;\mat{\mmul}{\opsp}\bigr)$,
$a\in\opsp$ and $t\in\R_+$ then
\[
\langle j^{\psi^\dagger}_t(a^*) u\evec{f},v\evec{g}\rangle=%
\langle u\evec{f},j^\psi_t(a) v\evec{g}\rangle\qquad
\forall\,u,v\in\ini,\ f,g\in\elltwo,
\]
by Theorem~\ref{thm:regcoc} and (\ref{eqn:viip}). Consequently,
$j^{\psi^\dagger}_t(a^*)\subseteq j^\psi_t(a)^*=:%
(j^\psi_t)^\dagger(a^*)$.
\end{rem}

\begin{rem}\label{rmk:hpee}
Let $Y$ be a vacuum-adapted quantum semimartingale on $\ffock$ which
satisfies the quantum stochastic differential equation
\begin{equation}\label{eqn:hpee1}
Y_0=\id_\ini\otimes\expn_0|_\eevecs,\qquad%
\rd Y_t=\rd\Lambda_F(t)\,Y_t\biggl(=\sum_{\alpha,\beta\in I}%
(F^\alpha_\beta\otimes\id_\fock)Y_t\std\Lambda^\beta_\alpha(t)\biggr),
\end{equation}
where $\smash[t]{F\in\bop{\ini\otimes\mmul}{}}$, \ie $Y$ is a
vacuum-adapted $\ini$~process such that
\begin{equation}\label{eqn:hpeeip}
\langle u\evec{f},Y_t v\evec{g}\rangle=%
\langle u,v\rangle+\int_0^t\langle u\evec{f},(E^{\widehat{f(s)}}F %
E_{\widehat{g(s)}}\otimes\id_\fock)Y_s v\evec{g}\rangle\std s
\end{equation}
for all $u$,~$v\in\ini$, $f$,~$g\in\elltwo$ and $t\in\R_+$; it follows
automatically that each $Y_t$ is bounded and $Y$ has locally uniformly
bounded norm. If $j_t(a):=(a\otimes\id_\fock)Y_t$ for all
$a\in\bop{\ini}{}$ and $t\in\R_+$ then it is readily verified
that~$j=j^\psi$, the vacuum cocycle with generator
\begin{equation}\label{eqn:hpgen}
\psi:\bop{\ini}{}\to\bop{\ini\otimes\mmul}{};\ %
a\mapsto(a\otimes\id_\mmul)F.
\end{equation}
Conversely, if $j$ is the vacuum cocycle with generator $\psi$ of the
form (\ref{eqn:hpgen}) then $Y=\bigl(j_t(\id_\ini)\bigr)_{t\in\R_+}$
is a vacuum-adapted quantum semimartingale which satisfies
(\ref{eqn:hpee1}). (If $\psi$ has the form (\ref{eqn:hpgen}) then
$\psi^{(m)}(a)=(a\otimes\id_{\mmul^{\otimes m}})\psi^{(m)}(\id_\ini)$
for all $m\in\Z_+$ and $a\in\bop{\ini}{}$. Hence
$j_t(a)=(a\otimes\id_\fock)j_t(\id_\ini)$ for all $t\in\R_+$, by
Theorem~\ref{thm:regcoc}, and (\ref{eqn:hpeeip}) holds.)

The other Hudson--Parthasarathy evolution equation,
\begin{equation}\label{eqn:hpee2}
\rd Z_t=Z_t\std\Lambda_G(t)\biggl(=\sum_{\alpha,\beta\in I}%
Z_t(G^\alpha_\beta\otimes\id_\fock)\std\Lambda^\beta_\alpha(t)\biggr),
\end{equation}
may, of course, be dealt with in a similar manner; if
$\smash[t]{G\in\bop{\ini\otimes\mmul}{}}$ then the corresponding
generator has the form $a\mapsto G(a\otimes\id_\mmul)$. The
observation that (\ref{eqn:flow}) generalises (\ref{eqn:hpee1}) and
(\ref{eqn:hpee2}) goes back at least to \cite[Proposition~5.2]{LiP98}.
\end{rem}

\section{Approximating walks}\label{sec:wlk}

\begin{defn}
If $n\in\Z_+$, $m\in\{0,\ldots,n\}$ and
$\bp\in\{0,\ldots,n-1\}^{m,\uparrow}$ then the normal $*$-homomorphism
$\vhm_\bp^n:\bop{\bbebe_{m)}}{}\to\bop{\bbebe_{n)}}{}$ is the
truncation of $\vhm_\bp$ described in Definition~\ref{def:vhm},
\ie $\vhm_\bp(a)=\vhm_\bp^n(a)\otimes\Delta^\perp_{[n}$ for all
$a\in\bop{\bbebe_{m)}}{}$.
\end{defn}

\begin{lem}\label{lem:modf}
Given $\phi\in\mmul\bopp\bigl(\opsp;\mat{\mmul}{\opsp}\bigr)$, define
$\modf{\phi}\in\mmul\bopp\bigl(\opsp;\mat{\mmul}{\opsp}\bigr)$ by
setting $\modf{\phi}(a):=\phi(a)-a\otimes\Delta^\perp$ for all
$a\in\opsp$. If $n\in\Z_+$ and $a\in\opsp$ then
\begin{equation}\label{eqn:modf}
\phi^{(n)}(a)=%
\sum_{m=0}^n\sum_{\bp\in\{0,\ldots,n-1\}^{m,\uparrow}}%
\vhm_\bp^n\bigl(\modf{\phi}^{(m)}(a)\bigr).
\end{equation}
\end{lem}
\begin{proof}
The cases $n=0$ and $n=1$ are immediately verified; suppose that
the identity (\ref{eqn:modf}) holds for some $n\geq1$. If
$x$,~$y\in\mmul$ and $a\in\opsp$ then
\begin{align*}
E^x\phi^{(n+1)}(a)E_y&=\phi^{(n)}(E^x\phi(a)E_y)\\
&=\sum_{m=0}^n\sum_{\bp\in\{0,\ldots,n-1\}^{m,\uparrow}}%
\vhm_\bp^n\bigl(\modf{\phi}^{(m)}(E^x\phi(a)E_y)\bigr)\\
&=\sum_{m=0}^n\sum_\bp%
\vhm_\bp^n\bigl(E^x\modf{\phi}^{(m+1)}(a)E_y+%
\langle x,\Delta^\perp y\rangle \modf{\phi}^{(m)}(a)\bigr)\\
&=\sum_{m=0}^n\sum_\bp E^x%
\Bigl(\vhm_{\bp\cup n}^{n+1}\bigl(\modf{\phi}^{(m+1)}(a)\bigr)+%
\vhm_\bp^{n+1}\bigl(\modf{\phi}^{(m)}(a)\bigr)\Bigr)E_y\\
&=%
E^x\Biggl(\sum_{m=0}^{n+1}\sum_{\bq\in\{0,\ldots,n\}^{m,\uparrow}}%
\vhm_\bq^{n+1}\bigl(\modf{\phi}^{(m)}(a)\bigr)\Biggr)E_y,
\end{align*}
where $\bp\cup n:=(p_1,\ldots,p_m,n)$ for all
$\bp=(p_1,\ldots,p_m)\in\{0,\ldots,n-1\}^{m,\uparrow}$; the result
follows by induction.
\end{proof}

\gap
The next proposition is the quantum-random-walk version of the
Wiener--It\^o decomposition given by Theorem~\ref{thm:regcoc}. Note
that each sum in (\ref{eqn:dwi}) has only finitely many non-zero terms
(as $\vint^m(\cdot;P_{(h)},h)_t=0$ if $m>t/h$).

\begin{prop}\label{prp:qrwwi}
If $h>0$ and $\phi\in\mmul\bopp\bigl(\opsp;\mat{\mmul}{\opsp}\bigr)$
then
\begin{align}\label{eqn:dwi}
J^{\phi,h}_t(a)&:=\sum_{n=0}^\infty\tfn{t\in{[n h,(n+1)h[}}%
D_h^*(\phi^{(n)}(a)\otimes\Delta^\perp_{[n})D_h\\
&=a\otimes\expn_0+%
\sum_{m=1}^\infty\vint^m(\scale_h(\modf{\phi})^{(m)}(a)%
\otimes\expn_\Vac;P_{(h)},h)_t\nonumber
\end{align}
for all $a\in\opsp$ and $t\in\R_+$. The mapping process $J^{\phi,h}$
is the \emph{vacuum-embedded random walk} with \emph{generator} $\phi$
and \emph{step size} $h$.
\end{prop}
\begin{proof}
This follows from Lemma~\ref{lem:modf} and
Propositions~\ref{prp:vdi} and \ref{prp:scale}, together with the
observation that $\{0,\ldots,n-1\}^{m,\uparrow}$ is empty whenever
$m>n$.
\end{proof}

\begin{rem}
Let $h>0$ and
$\phi\in\mmul\bopp\bigl(\opsp;\mat{\mmul}{\opsp}\bigr)$. In the same
manner as for Remark~\ref{rmk:cbcoc}, but using the inequality
(\ref{eqn:subonorm}) instead of (\ref{eqn:vinm}), it follows that
$\vint^m\bigl(\phi^{(m)}(\cdot)\otimes\expn_\Vac;P_{(h)},_h\bigr)_t %
E_{\evec{f}}\in\mmul\bopp\bigl(\opsp;\bop{\ini}{\ffock}\bigr)$,
with
\[
\|\vint^m\big(\phi^{(m)}(a)\otimes\expn_\Vac;P_{(h)},h)_t %
E_{\evec{f}}\|_\mmul \leq %
c_t^m \|\phi\|^m_\mmul \|\indf{[0,t[}\widehat{f}\|^m/\sqrt{m!},
\]
for all $m\geq1$, $t\in\R_+$ and $f\in\elltwo$, and the same inequality
holds if the $\mmul$ norm is replaced by the $\cb$ norm.
\end{rem}

\begin{rem}\label{rmk:invol3}
If $h>0$ and $\phi\in\mmul\bopp\bigl(\opsp;\mat{\mmul}{\opsp}\bigr)$
then Remark~\ref{rmk:walkadj} implies that
$J^{\phi^\dagger,h}_t(a^*)=J^{\phi,h}_t(a)^*=(J^{\phi,h}_t)^\dagger(a^*)$
for all $a\in\opsp$ and $t\in\R_+$.
\end{rem}

\begin{rem}
If
$\phi(a)=\left[\begin{smallmatrix}p&q\\r&s\end{smallmatrix}\right]$,
where $p\in\bop{\ini}{}$, $q\in\bop{\ini\otimes\mul}{\ini}$,
$r\in\bop{\ini}{\ini\otimes\mul}$ and $s\in\bop{\ini\otimes\mul}{}$, then
\[
\scale_h(\modf{\phi})(a)=\begin{bmatrix}(p-a)/h&q/\sqrt{h}\,\\
r/\sqrt{h}&s\end{bmatrix}.
\]
Thus the scaling employed in the following theorem is to be expected:
the components $p$, $q$, $r$ and $s$ correspond respectively to time,
annihilation, creation and gauge parts of the quantum stochastic
integral and the Poisson and Brownian increments are such that
$(\rd\Lambda)^2=\rd\Lambda$ and $(\rd B)^2=\rd t$, whence
$\rd\Lambda=\rd B/\sqrt{\rd t}=\rd t/\rd t$.
\end{rem}

\begin{thm}\label{thm:main}
Let $h_n>0$ and
$\phi_n$,~$\psi\in\mmul\bopp\bigl(\opsp;\mat{\mmul}{\opsp}\bigr)$ be
such that
\begin{equation}\label{eqn:cond}
h_n\to0\quad\mbox{and}\quad%
\lift{\scale_{h_n}(\modf{\phi_n})}{\mmul}\to%
\lift{\psi}{\mmul}\mbox{ strongly}
\end{equation}
 as $n\to\infty$. If $f\in\elltwo$ and $T\in\R_+$ then
\[
\lim_{n\to\infty}\sup_{t\in[0,T]}\|J^{\phi_n,h_n}_t(a)E_{\evec{f}}-%
j^\psi_t(a)E_{\evec{f}}\|=0\qquad\forall\,a\in\opsp;
\]
if, further,
$\lim\limits_{n\to\infty}\|\scale_{h_n}(\modf{\phi_n})-\psi\|_\mmul=0$
then
\begin{equation}\label{eqn:main2}
\lim_{n\to\infty}%
\sup_{t\in[0,T]}\|J^{\phi_n,h_n}_t(\cdot)E_{\evec{f}}-%
j^\psi_t(\cdot)E_{\evec{f}}\|_\mmul=0,
\end{equation}
and similarly with $\|\cdot\|_\mmul$ replaced by $\|\cdot\|_\cb$ if
$\phi_n$ and $\psi$ are completely bounded.
\end{thm}
\begin{proof}
For convenience, let $\zeta_n:=\scale_{h_n}(\modf{\phi_n})$. Then
$\{\|\zeta_n\|_\mmul:n\geq1\}$ is bounded (either by
Proposition~\ref{prp:liftcng} and the principle of uniform boundedness
or by Lemma~\ref{lem:cbequiv}) so, by Theorem~\ref{thm:regcoc},
Proposition~\ref{prp:qrwwi} and the inequalities (\ref{eqn:qinm}) and
(\ref{eqn:vicbnm}), it suffices to prove that, for all $m\geq1$ and
$f\in\elltwo$,
\begin{equation}\label{eqn:mth}
\|\bigl(\vint^m(\zeta_n^{(m)}(a)\otimes\expn_\Vac;P_{(h_n)},h_n)_t-%
\vint^m(\psi^{(m)}(a)\otimes\expn_\Vac)_t\bigr)E_{\evec{f}}\|\to0
\end{equation}
as $n\to\infty$ in the appropriate sense.

For this, note first that the left-hand side of (\ref{eqn:mth}) is
dominated by
\begin{multline*}
\|\bigl(\vint^m(\zeta_n^{(m)}(a)\otimes\expn_\Vac;P_{(h_n)},h_n)_t-%
\vint^m(\zeta_n^{(m)}(a)\otimes\expn_\Vac)_t\bigr)E_{\evec{f}}\|\\
+\|\vint^m\bigl((\zeta_n^{(m)}-%
\psi^{(m)})(a)\otimes\expn_\Vac\bigr)_t E_{\evec{f}}\|.
\end{multline*}
It follows from \cite[Propositions~3.2 and~3.3]{Blt07} that, if
$d_t:=2c_t^2=4\max\{t,1\}$,
\begin{align}\label{eqn:part1}
\|\bigl(\vint^m(\zeta_n^{(m)}(a)&\otimes\expn_\Vac;P_{(h)},h)_t-%
\vint^m(\zeta_n^{(m)}(a)\otimes\expn_\Vac)_t\bigr)u\evec{f}\|^2\\
&\leq d_t^m\int_{\Delta_m(t)\setminus\Delta_m^h(t)}%
\|(\zeta_n^{(m)}(a)\otimes\expn_0)\nnabla^{P_{(h)},m}_\bt %
u\evec{f}\|^2\std\bt\nonumber\\
&\qquad+d_t^m\|\evec{f}\|^2\sum_{n=1}^m%
\int_{\Delta_m(t)}\|\zeta_n^{(m)}(a)u\otimes g_n(h,\bt)\|^2\std\bt%
\nonumber
\end{align}
for all $h>0$ and $u\in\ini$, where
\[
g_n(h,\bt):=\widehat{P_{(h)}f}^{\otimes n-1}(t_1,\ldots,t_{n-1})%
\otimes(P_{(h)}-I)f(t_n)\otimes%
\widehat{f}^{\otimes m-n}(t_{n+1},\ldots,t_m).
\]
The first term on the right-hand side of (\ref{eqn:part1}) is
dominated by
\[
d_t^m\|\zeta_n\|_\mmul^{2m}\|a\|^2\|u\|^2%
\int_{\Delta_m(t)\setminus\Delta_m^h(t)}%
\|\widehat{P_{(h)}f}^{\otimes}(\bt)\|^2\std\bt
\]
and, as shown in Appendix~\ref{apx:estimate},
\[
\lim_{h\to0}\sup_{t\in[0,T]}%
\int_{\Delta_m(t)\setminus\Delta_m^h(t)}%
\|\widehat{P_{(h)}f}^{\otimes m}(\bt)\|^2\std\bt=0.
\]
Since
$\Delta^m(t)\subseteq\Delta^{n-1}(t)\times{[0,t[}\times\Delta^{m-n}(t)$,
the sum in (\ref{eqn:part1}) is dominated by
$\|\zeta_n\|_\mmul^{2m}\|a\|^2\|u\|^2 R_m$, where
\begin{multline*}
R_m:=%
\sum_{n=1}^m\frac{1}{(n-1)!}(t+\|P_{(h)}f\|^2)^{n-1}\|(P_{(h)}-I)f\|^2%
\frac{1}{(m-n)!}(t+\|f\|^2)^{m-n}\\
\leq \frac{1}{(m-1)!}2^{m-1}(t+\|f\|^2)^{m-1}\|(P_{(h)}-I)f\|^2.
\end{multline*}
Using (\ref{eqn:vicb}) (and the analogous result for the modified
integral) as necessary, it follows that (\ref{eqn:part1}) converges
to zero in the required manner.

Finally, if $X\in\bop{\ffock_{m)}}{}$ then (\ref{eqn:vinm}) implies
that
\[
\sup_{t\in[0,T]}\|\vint^m(X\otimes\expn_\Vac)_t E_{\evec{f}}\|\leq %
c_T^m \|\indf{[0,T[}\widehat{f}\|^m\|X\|/\sqrt{m!},
\]
and the result follows, by Lemma~\ref{lem:approxgen} and
(\ref{eqn:vicb}) (for the completely bounded case).
\end{proof}

\begin{rem}
When $\mul$ is infinite dimensional,
$\phi_n\matten\id_{\bop{\mmul}{}}\to0$ strongly if and only if
$\|\phi_n\|_\cb\to0$, by Lemma~\ref{lem:cbequiv}, so the stronger,
completely bounded version of (\ref{eqn:main2}) will hold in this
case.
\end{rem}

\begin{rem}
Proposition~\ref{prp:hsineq} shows that the `Hilbert--Schmidt'
condition used by Attal and Pautrat \cite[Remark on p.80]{AtP06} is
stronger than the convergence hypotheses required in
Theorem~\ref{thm:main}. They considered generators of the form
\[
\phi:\bop{\ini}{}\to\bop{\ini\otimes\mmul}{};\ %
a\mapsto(a\otimes\id_\mmul)L,
\]
where $L\in\bop{\ini\otimes\mmul}{}$ (\cf Remark~\ref{rmk:hpee}), so
$\lift{\phi}{\mmul}(b)=U^*(b\otimes\id_\mmul)U(L\otimes\id_\mmul)$ for
all $b\in\bop{\ini\otimes\mmul}{}$ (with
$U\in\bop{\ini\otimes\mmul\otimes\mmul}{}$ the unitary operator which
exchanges the last two tensor components). If $\{e_i:i\in I\}$ is an
orthonormal basis for $\mmul$ then (\ref{eqn:hsineq}) implies that
\begin{equation}\label{eqn:aphsineq}
\|\lift{\phi}{\mmul}\|^2\leq\|L\otimes\id_\mmul\|^2=\|L\|^2\leq%
\sum_{i,j\in I}\|L^i_j\|^2,
\end{equation}
where $L^i_j:=E^{e_i}L E_{e_j}$; the quantity on the right-hand side
of (\ref{eqn:aphsineq}) is that used by Attal and Pautrat for
establishing their results.
\end{rem}

\begin{thm}\label{thm:algmain}
Suppose $h_n>0$ and
$\phi_n$,~$\psi\in\mmul\bopp\bigl(\csa;\mat{\mmul}{\csa}\bigr)$ are
such that~(\ref{eqn:cond}) holds as $n\to\infty$, where the
$C^*$~algebra $\csa$ acts on $\ini$. If
\[
\phi_n=\phi_n^\dagger\quad\mbox{and}\quad%
\phi_n^{(m)}(ab)=\phi_n^{(m)}(a)\phi_n^{(m)}(b)\quad%
\forall\,a,b\in\csa,\ m,n\geq1
\]
(\ie $\phi_n$ is real and $\phi_n^{(m)}$ is multiplicative) then
$\psi=\psi^\dagger$ and the map $j^\psi_t$ is a $*$-homomorphism from
$\csa$ into $\bop{\ffock}{}$ for all $t\in\R_+$.
\end{thm}
\begin{proof}
The first claim holds by Proposition~\ref{prp:liftcng} and norm
continuity of the operator adjoint. Next, observe that the map
$a\mapsto J^{\phi_n,h_n}_t(a)$ is a $*$-homomorphism for all $n\geq1$
and $t\in\R_+$, which is therefore contractive, and so, by
Theorem~\ref{thm:main},
\[
\|j^\psi_t(a)\eta\|=%
\lim_{n\to\infty}\|J^{\phi_n,h_n}_t(a)\eta\|\leq%
\|a\|\,\|\eta\|\qquad\forall\,\eta\in\eevecs.
\]
Hence $j^\psi_t(a)\in\bop{\ffock}{}$ for all $a\in\csa$ and
$t\in\R_+$ (or, rather, $j^\psi_t(a)$ extends to such an operator,
but this distinction will be neglected). Remark~\ref{rmk:invol2} now
yields the identity $j^\psi_t(a^*)=j^\psi_t(a)^*$ for all
$t\in\R_+$ and $a\in\csa$; since
\begin{align*}
\langle\eta,j^\psi_t(ab)\zeta\rangle&=%
\lim_{n\to\infty}\langle\eta,J^{\phi_n,h_n}_t(ab)\zeta\rangle\\
&=\lim_{n\to\infty}\langle J^{\phi_n,h_n}_t(a^*)\eta,%
J^{\phi_n,h_n}_t(b)\zeta\rangle=%
\langle j^\psi_t(a^*)\eta,j^\psi_t(b)\zeta\rangle
\end{align*}
for all $\eta$,~$\zeta\in\eevecs$, $a$,~$b\in\csa$ and $t\in\R_+$,
the result follows.
\end{proof}

\begin{rem}\label{rmk:product}
If
\begin{enumerate}
\item[(a)] $\phi$ is a $*$-homomorphism from the $C^*$~algebra $\csa$
into $\csa\splten\bop{\mmul}{}$ or
\item[(b)] $\phi$ is a normal $*$-homomorphism from the von~Neumann
algebra $\csa$ into $\csa\uwkten\bop{\mmul}{}$
\end{enumerate}
then $\phi^{(m)}$ is multiplicative for all $m\geq1$, by
Proposition~\ref{prp:hom}. The hypotheses (a) and (b) are the
conditions ($\alpha_2$) and ($\beta_2$), respectively, of Lindsay and
Wills \cite[Corollary~4.2]{LiW03a}; in fact, the weaker conditions
($\alpha_1$) and ($\beta_1)$ given there (and below) each suffice to
establish multiplicativity. If the maps
$\phi\in\mmul\bopp\bigl(\csa;\mat{\mmul}{\csa}\bigr)$ and $\phi^{(m)}$
are both multiplicative for some $m\geq1$ and $\{e_i:i\in I\}$ is an
orthonormal basis for $\mmul$ then
\begin{align*}
E^x\phi^{(m+1)}(a)\phi^{(m+1)}(b)E_y&=\sum_{i\in I}%
\phi^{(m)}(E^x\phi(a)E_{e_i})\phi^{(m)}(E^{e_i}\phi(b)E_y)\\
&=\phi^{(m)}(E^x\phi(ab)E_y)
\end{align*}
for all $x$,~$y\in\mmul$ and $a$,~$b\in\csa$ provided either
\begin{enumerate}
\item[($\alpha_1$)] $\phi(c)E_z\in\csa\splten\bop{\C}{\mmul}$ for
all $c\in\csa$ and $z\in\mmul$ (in which case the series
$\sum_{i\in I}E_{e_i}E^{e_i}\phi(b)E_y$ is norm convergent) or
\item[($\beta_1$)] $\phi^{(m)}$ is strong-operator to weak-operator
continuous on bounded sets (\ie if the net $(c_\lambda)\subseteq\csa$
is norm bounded and converges in the strong operator topology to
$c\in\csa$ then $\phi^{(m)}(c_\lambda)\to\phi^{(m)}(c)$ in the weak
operator topology; if $\mmul$ is separable then `net' may be
replaced by `sequence').
\end{enumerate}
\end{rem}

\section{Semigroup dilation}\label{sec:exa}

\begin{rem}
Let the operator space $\opsp\subseteq\bop{\ini}{}$. If
$\phi\in\cbo{\opsp}{\mat{\hilb}{\opsp}}$ and
$R$,~$S\in\bop{\ini\otimes\hilb}{}$ then the map
\[
\psi:\opsp\to\mat{\hilb}{\opsp};\ a\mapsto R^*\phi(a)S
\]
is completely bounded, with
$\|\psi\|_\cb\leq\|R\|\,\|\phi\|_\cb\|S\|$, and, for any Hilbert space
$\hilc$, the map
$\lift{\psi}{\hilc}:\mat{\hilc}{\opsp}\to%
\mat{\hilb\otimes\hilc}{\opsp}$
has the form
\[
b\mapsto(R^*\otimes\id_\hilc)%
\bigl(\lift{\phi}{\hilc}(b)\bigr)(S\otimes\id_\hilc).
\]
If $\opsp$ is ultraweakly closed, so that
$\mat{\hilb}{\opsp}=\opsp\uwkten\bop{\hilb}{}$, and $\phi$ is
ultraweakly continuous then so are $\psi$ and
$\lift{\psi}{\hilc}=\psi\uwkten\id_{\bop{\hilc}{}}$.
\end{rem}

\begin{exmp}{[\Cf\cite[Theorem~4.2]{Sah08}.]}\label{exp:lind}
Let $\csa\subseteq\bop{\ini}{}$ be a unital $*$-algebra and suppose
either
\begin{enumerate}
\item[(i)] $\csa$ is a $C^*$~algebra, the map
$\pi:\csa\to\csa\splten\bop{\mul}{}$ is a $*$-homomorphism, $g\in\csa$
is self adjoint, $r\in\csa\splten\bop{\mul}{\C}$ and
$w\in\csa\splten\bop{\mul}{}$ is a co-isometry or
\item[(ii)] $\csa$ is a von~Neumann algebra, the map
$\pi:\csa\to\csa\uwkten\bop{\mul}{}$ is a normal $*$-homomorphism,
$g\in\csa$ is self adjoint, $r\in\csa\uwkten\bop{\mul}{\C}$ and
$w\in\csa\uwkten\bop{\mul}{}$ is a co-isometry.
\end{enumerate}
Let $\agb:=\csa\splten\bop{\mmul}{}$ in case~(i) and
$\agb:=\csa\uwkten\bop{\mmul}{}$ in case~(ii), noting that
$\agb\subseteq\mat{\mmul}{\csa}$. The $*$-homomorphism
\[
\phi_h:\csa\to\agb;\ %
a\mapsto V_h^*U_h^*\widehat{\pi}(a)U_h V_h\qquad(h>0)
\]
is normal in case~(ii), where
\[
\widehat{\pi}(a):=\begin{bmatrix}a&0\\0&\pi(a)\end{bmatrix},\ %
R:=\begin{bmatrix}0&-r\\r^*&0\end{bmatrix},\ %
U_h:=\e^{h^{1/2}R}\mbox{ and }%
V_h:=\begin{bmatrix}\e^{\I h g}&0\\0&w\e^{\I h g\otimes\id_\mul}%
\end{bmatrix}
\]
each belong to $\agb$. If
\begin{align}\label{eqn:homgen}
\psi(a)&:=%
B^*\widehat{\pi}(a)B+A^*\widehat{\pi}(a)C+C^*\widehat{\pi}(a)A\\
&\phantom{:}=\begin{bmatrix}
\I[a,g]-\half\{a,r r^*\}+r\pi(a)r^*&r\pi(a)w-a r w\\
w^*\pi(a)r^*-w^* r^* a&w^*\pi(a)w\end{bmatrix}\nonumber
\end{align}
for all $a\in\csa$, where $[\cdot,\cdot]$ is the commutator,
$\{\cdot,\cdot\}$ the anti-commutator,
\[
A:=\begin{bmatrix}\id_\ini&0\\0&0\end{bmatrix},\quad%
B:=\begin{bmatrix}0&0\\r^*&w\end{bmatrix}\quad\mbox{and}\quad %
C:=\begin{bmatrix}\I g-\half r r^*&-r w\\0&0\end{bmatrix},
\]
then $\|\scale_h(\modf{\phi_h})(a)-\psi(a)\|\to0$; note that
\[
\scale_h(\modf{\phi_h})(a)=\scale_h\bigl(\modf{\phi_h}(a)\bigr)=%
W_h^*\widehat{\pi}(a)W_h-a\otimes h^{-1}\Delta^\perp,
\]
where
\[
W_h:=U_h V_h(\id_\ini\otimes\Xi_h)=h^{-1/2}A+B+h^{1/2}C+O(h)
\]
as $h\to0$ (with respect to the norm topology). If $\hilc$ is a
Hilbert space and the unitary operator
$U\in\bop{\ini\otimes\mmul\otimes\hilc}%
{\ini\otimes\hilc\otimes\mmul}$ exchanges the last two tensor
components then
\begin{multline*}
\lift{\scale_h(\modf{\phi_h})}{\hilc}(b)\\
=(W_h\otimes\id_\hilc)^*\bigl(\lift{\widehat{\pi}}{\hilc}(b)\bigr)%
(W_h\otimes\id_\hilc)-U^*(b\otimes h^{-1}\Delta^\perp)U
\end{multline*}
for all $b\in\mat{\hilc}{\csa}$, so
$\lift{\scale_h(\modf{\phi_h})}{\hilc}\to\lift{\psi}{\hilc}$ strongly
as $h\to0$; Theorem~\ref{thm:algmain} implies that $\psi$ is the
generator of a $*$-homomorphic vacuum cocycle. Furthermore,
\[
T_t:\csa\to\csa;\ a\mapsto E^{\evec{0}}j^\psi_t(a) E_{\evec{0}}%
\qquad(t\in\R_+)
\]
defines a semigroup of completely positive contractions (a
\emph{CPC semigroup}) with bounded generator
\[
\lind:a\mapsto\I[a,g]-\half\{a,r r^*\}+r \pi(a) r^*,
\]
which is ultraweakly continuous (and so each $T_t$ is also) in
case~(ii).
\end{exmp}

\begin{exmp}
Conversely, suppose $(T_t)_{t\in\R_+}$ is a norm-continuous CPC
semigroup of ultraweakly continuous operators on the von~Neumann
algebra~$\csa$. There exists an ultraweakly continuous operator
$\lind\in\bop{\csa}{}$ such that $T_t=\exp(t\lind)$ for all
$t\in\R_+$; moreover, there exist a Hilbert space $\mul$, a normal
$*$-homomorphism $\pi:\csa\to\csa\uwkten\bop{\mul}{}$, a self-adjoint
operator $g\in\csa$ and an operator $r\in\csa\uwkten\bop{\mul}{\C}$
such that
\begin{equation}\label{eqn:lind}
\lind(a)=%
\I[a,g]-\half\{a,r r^*\}+r \pi(a) r^*\qquad\forall\,a\in\csa.
\end{equation}
(See \cite[Theorem~6.9]{Lin05}; this is the Goswami--Sinha version
\cite[Theorem~4.0.1]{GoS99} of the
Gorini--Kossakowski--Sudarshan/Lindblad decomposition established by
Christensen and Evans.) Example~\ref{exp:lind} yields a
$*$-homomorphic vacuum cocycle $j^\psi$, where $\psi$ is defined as in
(\ref{eqn:homgen}) (with $w\in\csa\uwkten\bop{\mul}{}$ an arbitrary
co-isometry), such that
\[
T_t(a)=E^{\evec{0}}j^\psi_t(a)E_{\evec{0}}\qquad%
\forall\,t\in\R_+,\ a\in\csa.
\]
\end{exmp}

\begin{exmp}
Suppose alternatively that $\lind\in\bop{\csa}{}$ is the generator of
a CPC semigroup on the unital separable $C^*$~algebra
$\csa\subseteq\bop{\ini}{}$. There exist a separable Hilbert space
$\mul$, a $*$-homomorphism $\pi:\csa\to\bop{\ini\otimes\mul}{}$ such
that
\[
\pi(a)E_x\in\csa\splten\bop{\C}{\mul}\qquad%
\forall\,a\in\csa,\ x\in\mul
\]
(whence $\pi(\csa)\subseteq\mat{\mul}{\csa}$) and a $\pi$-derivation
$\delta:\csa\to\csa\splten\bop{\C}{\mul}$ (\ie a linear map
satisfying $\delta(ab)=\delta(a)b+\pi(a)\delta(b)$ for all
$a$,~$b\in\csa$) such that
\[
\delta(a)^*\delta(b)=%
\lind(a^* b)-\lind(a)^* b-a^*\lind(b)\qquad\forall\,a,b\in\csa.
\]
(See \cite[Proof of Theorem~6.7]{Lin05}.) Working \textit{via} the
universal representation \cite[{\S}10.1]{KaR97}, or otherwise, it may
be shown there exists a normal $*$-homomorphism
$\widetilde{\pi}:A''\to\csa''\uwkten\bop{\mul}{}$ such that
$\widetilde{\pi}|_\csa=\pi$. The Christensen--Evans result on twisted
derivations \cite[Theorem~2.1]{ChE79} yields
$r\in\csa''\uwkten\bop{\mul}{\C}$ and $g=g^*\in\csa''$ such that the
maps
\[
\widetilde{\delta}:\csa''\to\csa''\uwkten\bop{\C}{\mul};\ %
a\mapsto\widetilde{\pi}(a)r^*-r^*a
\]
and
\[
\widetilde{\lind}:\csa''\to\csa'';\ %
a\mapsto\I[a,g]-\half\{a,r r^*\}+r\widetilde{\pi}(a)r^*
\]
satisfy $\widetilde{\delta}|_\csa=\delta$ and
$\widetilde{\lind}|_\csa=\lind$ (\ie (\ref{eqn:lind}) holds). Let
$w\in\csa\splten\bop{\mul}{}$ be a co-isometry and define
\[
\widetilde{\psi}:\csa''\to\mat{\mmul}{\csa''};\ a\mapsto%
\begin{bmatrix}\widetilde{\lind}(a)&\widetilde{\delta}(a^*)^*w\\
w^*\widetilde{\delta}(a)&w^*\widetilde{\pi}(a)w\end{bmatrix};
\]
Example~\ref{exp:lind} shows that $j^{\widetilde{\psi}}$ is a
$*$-homomorphic vacuum cocycle and therefore so is its restriction to
$\csa$, which has generator
\[
\psi:=\widetilde{\psi}|_{\csa}=%
a\mapsto\begin{bmatrix}\lind(a)&\delta(a^*)^* w\\
w^*\delta(a)&w^*\pi(a)w\end{bmatrix}\in\mat{\mmul}{\csa},
\]
and $\exp(t\lind)(a)=E^{\evec{0}}j^\psi_t(a)E_{\evec{0}}$ for all
$t\in\R_+$ and $a\in\csa$, as above.
\end{exmp}

\section{Switching adaptedness}\label{sec:ssu}

\begin{notation}
For all $t\in\R_+$, let $\fock_{t)}$ and $\fock_{[t}$ denote Boson
Fock space over $L^2({[0,t[};\mul)$ and $L^2({[t,\infty[};\mul)$,
respectively, and let $\ffock_{t)}:=\ini\otimes\fock_{t)}$. Note that
$\fock_{t)}$ and $\fock_{[t}$ are subspaces of $\fock$.

Let $U_t:\ffock_{t)}\otimes\fock_{[t}\to\ffock$ be the unitary
operator with
$U_t^{-1}u\evec{f}=u\evec{f|_{[0,t[}}\otimes\evec{f|_{[t,\infty[}}$
for all $u\in\ini$ and $f\in\elltwo$; recall that an $\ini$~process
$(X_t)_{t\in\R_+}$ is \emph{identity adapted} if
\[
X_t u\evec{\indf{[0,t[}f}\in\ffock_{t)}\quad\mbox{and}\quad%
X_t u\evec{f}=%
U_t\Bigl(\bigl(X_t u\evec{\indf{[0,t[}f}\bigr)\otimes%
\evec{f|_{[t,\infty[}}\Bigr)
\]
for all $u\in\ini$, $f\in\elltwo$ and $t\in\R_+$.

Let
$U_t':=U_t(\widetilde{P}_{t)}\otimes P_{[t})\in%
\bop{\ffock\otimes\fock}{\ffock}$, where $\widetilde{P}_{t)}$ is the
orthogonal projection onto $\ffock_{t)}$ \etc, so
$U_t'(u\evec{f}\otimes\evec{g})=%
u\evec{\indf{[0,t[}f+\indf{[t,\infty[}g}$
for all $u\in\ini$,$f$,~$g\in\elltwo$ and $t\in\R_+$.
\end{notation}

\begin{thm}{\rm\bf[Lindsay \& Wills]}
If the operator space $\opsp\subseteq\bop{\ini}{}$ and
$\theta\in\mmul\bopp\bigl(\opsp;\mat{\mmul}{\opsp}\bigr)$ then there
exists a unique weakly regular mapping process $k^\theta$, the
\emph{QS cocycle generated by $\theta$}, such that, for all
$a\in\opsp$, the $\ini$~process $t\mapsto k^\theta_t(a)$ is identity
adapted and
\begin{equation}\label{eqn:ucocycle}
\langle u\evec{f},(k^\theta_t(a)-a\otimes\id_\fock)v\evec{g}\rangle=%
\int_0^t\langle u\evec{f},k^\theta_s(E^{\widehat{f(s)}}\theta(a)%
E_{\widehat{g(s)}}) v\evec{g}\rangle\std s
\end{equation}
for all $u$,~$v\in\ini$, $f$,~$g\in\elltwo$ and $t\in\R_+$.
\end{thm}
\begin{proof}
See \cite{LiW01} (for the case when $\opsp$ is a unital $C^*$~algebra;
as observed in \cite[Remark~3.2(i)]{LiW03a}, the extension to an
operator space is immediate).
\end{proof}

\begin{thm}\label{thm:switch}
Let
$\theta$,~$\psi\in\mmul\bopp\bigl(\opsp;\mat{\mmul}{\opsp}\bigr)$. If
$k^\theta$ is the QS cocycle with generator $\theta$ then $j$ is the
vacuum cocycle with generator $a\mapsto\theta(a)+a\otimes\Delta$,
where
\begin{equation}\label{eqn:jdef}
j_t(a) u\evec{f}:=k^\theta_t(a) u\evec{\indf{[0,t[}f}\qquad%
\forall\,u\in\ini,\ f\in\elltwo
\end{equation}
and $\Delta\in\bop{\mmul}{}$ is the orthogonal projection onto $\mul$.
Conversely, if $j^\psi$ is the vacuum cocycle with generator $\psi$
then $k$ is the QS cocycle with generator
$a\mapsto\psi(a)-a\otimes\Delta$, where
\[
k_t(a) u\evec{f}:=U_t\Bigl(\bigl(j^\psi_t(a)%
u\evec{\indf{[0,t[}f}\bigr)\otimes\evec{f|_{[t,\infty[}}\Bigr)%
\quad\forall\,u\in\ini,\ f\in\elltwo.
\]
\end{thm}
\begin{proof}
Note that first that if $a\in\opsp$ and $t\in\R_+$ then
$j_t(a) u\evec{\indf{[0,t[}f}\in\ffock_{t)}$ for all $u\in\ini$ and
$f\in\elltwo$, where $j$ is any (vacuum or QS) cocycle. If $k^\theta$
is the QS cocycle with generator $\theta$ and $j$ is defined by
(\ref{eqn:jdef}) then, by (\ref{eqn:ucocycle}),
\begin{align*}
\langle u\evec{f},j_t(a) v\evec{g}\rangle&=%
\langle u,a v\rangle\e^{\int_0^t\langle f(s),g(s)\rangle\std s}\\
&\ +\int_0^t\langle u\evec{f},k^\theta_s(E^{\widehat{f(s)}}\theta(a)%
E_{\widehat{g(s)}}) v\evec{g}\rangle\std s%
\,\e^{{-}\int_t^\infty\langle f(s),g(s)\rangle\std s}.
\end{align*}
Differentiation shows that
\[
\frac{\rd}{\rd t}\langle u\evec{f},j_t(a) v\evec{g}\rangle=%
\bigl\langle u\evec{f},j_t\bigl(E^{\widehat{f(t)}}\theta(a)%
E_{\widehat{g(t)}}+\langle f(t),g(t)\rangle a\bigr) %
v\evec{g}\bigr\rangle,
\]
which gives the first claim. The second may be established similarly.
\end{proof}

\begin{defn}
Let $h>0$, $\phi\in\mmul\bopp\bigl(\opsp;\mat{\mmul}{\opsp}\bigr)$ and
define
\[
K^{\phi,h}_t(a):=\sum_{n=0}^\infty\tfn{t\in[n h,(n+1)h[}%
D_h^*(\phi^{(n)}(a)\otimes\id_{\bebe_{[n}})D_h
\]
for all $a\in\opsp$ and $t\in\R_+$. The mapping process $K^{\phi,h}$
is the \emph{identity-embedded random walk} with \emph{generator}
$\phi$ and \emph{step size} $h$.
\end{defn}

\begin{notation}
If $\phi\in\kbo{\hilc}{\opsp}{\mat{\hilc}{\opsp}}$ then
$\mmodf{\phi}\in\kbo{\hilc}{\opsp}{\mat{\hilc}{\opsp}}$ is defined by
setting $\mmodf{\phi}(a):=\phi(a)-a\otimes\id_\hilc$ for all
$a\in\opsp$. This is the appropriate perturbation to a generator in
the identity-adapted situation; since
\[
\modf{\phi}(a)-a\otimes\Delta=%
\phi(a)-a\otimes\id_\mmul=\mmodf{\phi}(a),
\]
this perturbation is as expected \cite{Blt04}. The identity in
Lemma~\ref{lem:modf} is replaced by the following:
\[
\phi^{(n)}(a)=\sum_{m=0}^n\sum_{\bp\in\{0,\ldots,n-1\}^{m,\uparrow}}%
\uhm_\bp^n\bigl(\mmodf{\phi}^{(m)}(a)\bigr)\qquad%
\forall\,a\in\opsp,\ n\in\Z_+,
\]
where $\uhm_\bp^n(b)\otimes\id_{\bebe_{[n}}=\uhm_\bp(b)$ for all
$b\in\bop{\bbebe_{m)}}{}$ and
$\uhm_\bp:\bop{\bbebe_{m)}}{}\to\bop{\bbebe}{}$ is the normal
$*$-homomorphism such that
\[
X\otimes B_1\otimes\cdots\otimes B_m\mapsto X\otimes%
\id_{\bebe_{[0,p_1)}}\otimes B_1\otimes\id_{\bebe_{[p_1+1,p_2)}}%
\otimes\cdots\otimes B_m\otimes\id_{\bebe_{[p_m+1}},
\]
with $\bebe_{[p,q)}:=\bigotimes_{l=p}^{q-1}\mmul_{(l)}$ and
$\uhm_\emptyset(X)=X\otimes\id_\bebe$.
\end{notation}

\begin{thm}\label{thm:smain}
Suppose $h_n>0$ and
$\phi_n$,~$\theta\in\mmul\bopp\bigl(\opsp;\mat{\mmul}{\opsp}\bigr)$
are such that
\[
h_n\to0\quad\mbox{and}\quad%
\lift{\scale_{h_n}\bigl(\mmodf{\phi_n}\bigr)}{\mmul}%
\to\lift{\theta}{\mmul}\mbox{ strongly}
\]
as $n\to\infty$. If $f\in\elltwo$ and $T\in\R_+$ then
\[
\lim_{n\to\infty}\sup_{t\in[0,T]}%
\|K^{\phi_n,h_n}_t(a)E_{\evec{f}}-k^\theta_t(a)E_{\evec{f}}\|=0%
\qquad\forall\,a\in\opsp,
\]
where $k^\theta$ is the QS cocycle generated by $\theta$; if
$\lim\limits_{n\to\infty}%
\|\scale_{h_n}\bigl(\mmodf{\phi_n}\bigr)-\theta\|_\mmul=0$
then
\[
\lim_{n\to\infty}\sup_{t\in[0,T]}%
\|K^{\phi_n,h_n}_t(\cdot)E_{\evec{f}}-%
k^\theta_t(\cdot)E_{\evec{f}}\|_\mmul=0,
\]
and similarly with $\|\cdot\|_\mmul$ replaced by $\|\cdot\|_\cb$ if
$\phi_n$ and $\theta$ are completely bounded.
\end{thm}
\begin{proof}
Let $a\in\opsp$. Theorem~\ref{thm:switch} implies that
\[
k^\theta_t(a)E_{\evec{f}}=%
U_t' E_{\evec{f}} j^\psi_t(a)E_{\evec{f}},
\]
where $j^\psi$ is the vacuum cocycle with generator
$\psi:b\mapsto\theta(b)+b\otimes\Delta$, and if $t\in{[n h,(n+1)h[}$
then
\[
K^{\phi,h}_t(a)E_{\evec{f}}=%
U_{n h}' E_{D_h^*D_h\evec{f}} J^{\phi,h}_t(a)E_{\evec{f}},
\]
where $D_h$ acts on $\C\otimes\fock=\fock$. Hence
\begin{align*}
\|K^{\phi,h}_t(a)&E_{\evec{f}}-k^\theta_t(a)E_{\evec{f}}\|\\
&\leq\|(U_{n h}'-U_t')%
E_{D_h^*D_h\evec{f}} J^{\phi,h}_t(a)E_{\evec{f}}\|\\
&\quad+\|(D_h^*D_h-\id_\fock)\evec{f}\|\,%
\|J^{\phi,h}_t(a)E_{\evec{f}}\|\\
&\quad+\|\evec{f}\|\,%
\|J^{\phi,h}_t(a)E_{\evec{f}}-j^\psi_t(a)E_{\evec{f}}\|\\
&\leq\bigl(\|(P_{[n h}-P_{[t})D_h^*D_h\evec{f}\|+%
\|(D_h^*D_h-I_\fock)\evec{f}\|\bigr)\|J^{\phi,h}_t(a)E_{\evec{f}}\|\\
&\quad+\|\evec{f}\|\,%
\|J^{\phi,h}_t(a)E_{\evec{f}}-j^\psi_t(a)E_{\evec{f}}\|,
\end{align*}
where the second inequality holds because
$U_{n h}'E_\eta=U_t' E_{P_{[n h}\eta}$ for all $\eta\in\ffock$; the
same conclusion holds with the norm of the evaluation at $a$ is
replaced by $\|\cdot\|_\mmul$ or $\|\cdot\|_\cb$.
Theorem~\ref{thm:main} and the uniform strong continuity of
$t\mapsto P_{[t}$ (on $\evecs$, so everywhere) now give the result,
since if $a\in\opsp$ then
\[
(\scale_h(\modf{\phi})-\psi)(a)=%
\scale_h(\phi)(a)-\scale_h(a\otimes\Delta^\perp)-%
\theta(a)-a\otimes\Delta=%
(\scale_h\bigl(\mmodf{\phi}\bigr)-\theta)(a).\qedhere
\]
\end{proof}

\begin{exmp}
Let $\ini=\mul=\C$, $\opsp=\bop{\C}{}=\C I_\C$ and note the
isomorphism
\[
\bopp\bigl(\opsp;\opsp\matten\bop{\mmul}{}\bigr)\to\bop{\C^2}{};\ %
\phi\mapsto\phi(\id_\C).
\]
Fix $c\in\R$ and define $\phi$ and $\theta$ by setting
\[
\phi(\id_\C):=\left[\begin{array}{cc}
1&h^{1/2}\\
h^{1/2}&1+c
\end{array}\right]\quad(h>0)\qquad\mbox{and}\qquad%
\theta(\id_\C):=\left[\begin{array}{cc}
0&1\\
1&c
\end{array}\right].
\]
By Theorem~\ref{thm:smain}, 
$\sup_{t\in[0,T]}\|K^{\phi,h}_t(\id_\C)E_{\evec{f}}-%
k^\theta_t(\id_\C)E_{\evec{f}}\|\to0$
as $h\to0$, for all $f\in L^2(\R_+)$ and $T\in\R_+$, where
\[
\langle\evec{f},(k^\theta_t(\id_\C)-\id_\fock)\evec{g}\rangle=%
\int_0^t\langle\evec{f},%
k^\theta_s(\id_\C)\evec{g}\rangle%
\bigl(\overline{f(s)}+g(s)+c\overline{f(s)}g(s)\bigr)\std s
\]
for all $f$,~$g\in L^2(\R_+)$. Hence
$X=\bigl(k^\theta_t(\id_\C)\bigr)_{t\in\R_+}$ satisfies the quantum
stochastic differential equation
\[
X_0=\id_\evecs,\qquad%
\rd X_t=X_t(\rd A_t+\rd A^\dagger_t+c\std\Lambda_t)
\]
and therefore $X$ is unitarily equivalent to (multiplication by) the
Dol\'{e}ans-Dade exponential of the compensated Poisson process with
jump size $c$ and intensity $c^{-2}$ (if $c\neq0$) or of standard
Brownian motion (if $c=0$) \cite[p.74]{Mey95}. From this observation and
Yor's product formula \cite[Theorem~II.38]{Pro05} it follows that
$\evecs$ is stable under the action of $X$, with
\[
X_t\evec{f}=%
\e^{\int_0^t f(s)\std s}\evecc\bigl(f+\indf{[0,t[}(1+c f)\bigr)%
\qquad\forall\,f\in L^2(\R_+),\ t\in\R_+.
\]
Hence if $t\in\R_+$ and $m\in\Z_+$ then
\[
\langle\evec{0},X_t^m\evec{0}\rangle=%
\left\{\begin{array}{ll}
\exp\bigl(c^{-2}\bigl((1+c)^m-1-m c\bigr)t\bigr)&%
\mbox{if }c\neq0,\\[1ex]
\exp(m(m-1)t/2)&\mbox{if }c=0.
\end{array}\right.
\]

The operator process $Y=\bigl(K^{\phi,h}_t(\id_\C)\bigr)_{t\in\R_+}$
is commutative for all $h>0$, as
$\phi^{(n)}(\id_\C)=\phi(\id_\C)^{\otimes n}$ for all $n\in\Z_+$, and
\begin{equation}\label{eqn:walk1}
\langle\evec{0},Y_{n h}^m\evec{0}\rangle=%
\langle\vac,\phi(\id_\C)^m\vac\rangle^n\qquad\forall\,m,n\in\Z_+.
\end{equation}
The identity
$(\phi(\id_\C)-\id_{\C^2})^2=c(\phi(\id_\C)-\id_{\C^2})+h\id_{\C^2}$
may be used to show that
\begin{equation}\label{eqn:walk2}
\langle\vac,\phi(\id_\C)^m\vac\rangle=%
\Bigl(1-\frac{2h}{c+d}\Bigr)^m\frac{c+d}{2d}+%
\Bigl(1-\frac{2h}{c-d}\Bigr)^m\frac{d-c}{2d}%
\qquad\forall\,m\in\Z_+,
\end{equation}
where $d:=(c^2+4h)^{1/2}$. Thus $Y_{n h}$ is distributed in the vacuum
state as the product of $n$ independent copies of the random variable
$Z$, where
\[
\Pr\bigr(Z=(c+d-2h)/(c+d)\bigr)=p\quad\mbox{and}\quad%
\Pr\bigr(Z=(c-d-2h)/(c-d)\bigr)=1-p,
\]
with $p:=(c+d)/(2d)$. Hence, in this state, $(Y_{n h})_{n\geq0}$ is a
random walk with possible positions
\[
\Bigl(1-\frac{2h}{c+d}\Bigr)^j%
\Bigl(1-\frac{2h}{c-d}\Bigr)^{n-j}=%
\Bigl(\frac{c-d+2}{c+d+2}\Bigr)^j\Bigl(\frac{c-d-2h}{c-d}\Bigr)^n%
\qquad(j=0,\ldots,n)
\]
upon taking step $n$. (These are distinct unless $c=-2$; in this case,
\[
\Pr\bigl(Y_{n h}=\pm(1+h)^{n/2}\bigr)=1/2\qquad\forall\,n\geq1
\]
and the stochastic exponential which corresponds to $X$ takes the form
$t\mapsto\exp(t/2)(-1)^{N_t}$, where $N_t$ is the number of jumps
performed by the Poisson process up to time $t$.) After some working,
(\ref{eqn:walk1}) and (\ref{eqn:walk2}) give that
\[
\lim_{h\to0}\langle\evec{0},Y_t^m\evec{0}\rangle=%
\langle\evec{0},X_t^m\evec{0}\rangle%
\qquad\forall\,m\in\Z_+,\ t\in\R_+,
\]
\ie in the vacuum state, the moments of $Y_t$ converge to those of
$X_t$. However, the operators $X_s$ and $Y_t$ do not commute (if
$s>0$), so these processes are not simultaneously realisable: taking
$T$ to be sufficiently large, $s\in{[n h,(n+1)h[}$ and
$t\in{[p h,(p+1)h[}$,
\begin{multline*}
\langle k^\theta_s(\id_\C)\evec{\indf{[0,T[}},%
K^{\phi,h}_t(\id_\C)\evec{0}\rangle\\
=\left\{\begin{array}{ll}
\e^t(1+(2+c)h)^p&\mbox{if }n\geq p,\\[1ex]
\e^t(1+(2+c)h)^n\bigl(1+h+(1+c)(t-n h)\bigr)(1+h)^{p-n-1}&\mbox{if }n<p
\end{array}\right.
\end{multline*}
whereas
\begin{multline*}
\langle K^{\phi,h}_t(\id_\C)\evec{\indf{[0,T[}},%
k^\theta_s(\id_\C)\evec{0}\rangle\\
=\left\{\begin{array}{ll}
(1+(3+c)h)^n(1+h)^{n-p}(1+t-n h)&\mbox{if }n\geq p,\\[1ex]
(1+(3+c)h)^n\bigl(1+h+(2+c)(t-n h)\bigr)(1+h)^{p-n-1}&\mbox{if }n<p.
\end{array}\right.
\end{multline*}
\end{exmp}

\section*{Acknowledgements}
Part of this research was undertaken while the author was an Embark
Postdoctoral Fellow at University College Cork, funded by the Irish
Research Council for Science, Engineering and Technology. The
influence of Martin Lindsay and Steve Wills, through their work
\cite{LiW00,LiW01,LiW03a} and in many other ways, should be apparent
to the reader; the author is very happy to declare his debt to
them. Conversations with Lingaraj Sahu (during his visit to UCC and
elsewhere) and Adam Skalski are gratefully acknowledged. Robin Hudson
helpfully drew the author's attention to \cite{CGH77}.

\appendix

\section{An estimate}\label{apx:estimate}

Let $m\geq1$, $h>0$ and $f\in\elltwo$. Since $P_{(h)} f$ is
subordinate to the partition $\{0<h<2h<\cdots\}$, if
$g_h(p):=h\|\widehat{P_{(h)} f}(p h)\|^2\leq\|%
\indf{[p h,(p+1)h[}\widehat{f}\|^2$
then
\[
\int_{\Delta_m(n h)}%
\|\widehat{P_{(h)}f}^{\otimes m}(\bt)\|^2\std\bt=%
\frac{1}{m!}%
\sum_{\bp\in\{0,\ldots,n-1\}^m}g_h(p_1)\cdots g_h(p_m)
\]
and
\[
\int_{\Delta_m^h(n h)}%
\|\widehat{P_{(h)}f}^{\otimes m}(\bt)\|^2\std\bt=%
\frac{1}{m!}%
\sum_{\substack{\bp\in\{0,\ldots,n-1\}^m\\p_1,\ldots,p_m\text{ distinct}}}%
g_h(p_1)\cdots g_h(p_m).
\]
Let $\{0,\ldots,n-1\}^{m,\neq}$ denote the set of $m$-tuples of
distinct elements of $\{0,\ldots,n-1\}$ and
$\{0,\ldots,n-1\}^{m,=}:=%
\{0,\ldots,n-1\}^m\setminus\{0,\ldots,n-1\}^{m,\neq}$ the set of
$m$-tuples which contain at least one pair. If $t\in{[n h,(n+1)h[}$
then
\begin{align*}
I_m&:=m!\int_{\Delta_m(t)\setminus\Delta_m^h(t)}%
\|\widehat{P_{(h)}f}^{\otimes m}(\bt)\|^2\std \bt\\
&\leq\sum_{\bp\in\{0,\ldots,n\}^m\setminus\{0,\ldots,n-1\}^{m,\neq}}%
g_h(p_1)\cdots g_h(p_m)\\
&\leq m g_h(n)\|\indf{[0,(n+1)h[}\widehat{f}\|^{2(m-1)}+%
\sum_{\bp\in\{0,\ldots,n-1\}^{m,=}}g_h(p_1)\cdots g_h(p_m)\\
&\leq m g_h(n)\|\indf{[0,(n+1)h[}\widehat{f}\|^{2(m-1)}+%
\sum_{p=0}^{n-1} g_h(p)^2\|\indf{[0,n h[}\widehat{f}\|^{2(m-2)}.%
\end{align*}
Furthermore, $\sum_{p=0}^{n-1}g_h(p)^2$ is dominated by
\[
\sum_{p=0}^{n-1}\|\indf{[p h,(p+1)h[}\widehat{f}\|^4\leq%
\sup_{0\leq p\leq n-1}\|\indf{[p h,(p+1)h[}\widehat{f}\|^2%
\|\indf{[0,n h[}\widehat{f}\|^2,
\]
hence
\[
I_m\leq(m+1)\bigl(h+\sup_{p\geq0}\|\indf{[p h,(p+1)h[}f\|^2\bigr)%
(t+h+\|f\|^2)^{m-1}.
\]

\end{document}